\documentclass{article}
\usepackage{latexsym}
\usepackage{amsmath}
\usepackage{amscd}
\usepackage{amsthm}
\usepackage{amssymb}
\usepackage{color,enumitem,graphicx}
\usepackage{hyperref}

\usepackage{indentfirst}

\newtheorem{teorema}{Theorem}[section]
\newtheorem{proposicao}[teorema]{Proposition}
\newtheorem{corolario}[teorema]{Corollary}
\newtheorem{lema}[teorema]{Lemma}
\newtheorem{obs}[teorema]{Remark}

\usepackage[left=2.9cm,right=2.9cm,top=2.8cm,bottom=2.8cm]{geometry}

\newcommand{\fim}{\hfill $\Box$}

\def\fd{\hfill{\vbox to 7pt{\hbox to 7pt{\vrule height 7pt width 7pt}}}}

\def\N{\widetilde N}

\def\:{\colon}

\def\ds{\displaystyle}

   \def\NN{\mathbb{N}}
   \def\R{\mathbb{R}}
   \def\loc{\mathop{\rm loc}\nolimits}
   \def\D{{\nabla}}
   \def\lo{\mathop{\longrightarrow}}
   \def\cH{{\cal H}}
   
   \def\cB{{\cal B}}
   \def\cM{{\mathcal{M}}}

   \newcommand{\beq}{\begin{equation}}
   \newcommand{\eeq}{\end{equation}}

\title{Multiple positive bound state solutions of a critical  Choquard equation}

\author{Claudianor O. Alves
\\{\footnotesize Unidade Acad\^emica de Matem\'atica
, Universidade Federal de Campina Grande,}
\\ {\footnotesize 58429-970, Campina Grande - PB, Brazil.}
\\{\footnotesize E-mail address: coalves@mat.ufcg.edu.br}
\\
\\ Giovany M. Figueiredo
\\{\footnotesize  Departamento de Matem\'atica, Universidade de Brasilia - UNB}
\\{\footnotesize CEP: 70910-900, Bras\'ilia-DF, Brazil.}
\\{\footnotesize E-mail address: giovany@unb.br}
\\and\\
\\ Riccardo Molle
\\{\footnotesize  Dipartimento di Matematica, Universit\`a di Roma ``Tor Vergata"}
\\{\footnotesize CEP: 00133, Roma, Italia.}
\\{\footnotesize E-mail address: molle@mat.uniroma2.it}
}

\begin{document}

\pretolerance10000

\maketitle

\begin{abstract}
In this paper we consider the problem
 $$
 \left\{
\begin{array}{rcl}
-\Delta
u+V_{\lambda}(x)u=(I_{\mu}*|u|^{2^{*}_{\mu}})|u|^{2^{*}_{\mu}-2}u \ \
\mbox{in} \ \ \mathbb{R}^{N},\\ 
u>0 \ \ \mbox{in} \ \ \mathbb{R}^{N},
\end{array}
\right.\leqno{(P_{\lambda})}
$$
where $V_{\lambda}=\lambda+V_{0}$ with $\lambda \geq 0$,  $V_0\in
L^{N/2}(\R^N)$, $I_{\mu}=\frac{1}{|x|^\mu}$ is the Riesz potential
with $0<\mu<\min\{N,4\}$ and $2^{*}_{\mu}=\frac{2N-\mu}{N-2}$ with
$N\geq 3$. 
Under some smallness assumption on $V_0$ and $\lambda$ we prove the
existence of two positive solutions of $(P_\lambda)$. 
In order to prove the main result,  we used variational methods
combined with degree theory.    
\end{abstract}
\noindent{{\bf Key Words}. Choquard equation,  Variational methods,
  Critical exponents} 
\newline
\noindent{{\bf 2020 AMS Classification.} 81Q05, 35A15, 35B33} \vskip 0.4cm

\numberwithin{equation}{section}
\bibliographystyle{plain} 
\maketitle

\section{Introduction}

In this paper we will focus our attention on the existence of positive
solutions for the following class of Choquard equation  

$$
 \left\{
\begin{array}{rcl}
-\Delta u+(\lambda+V_0(x))u =
(I_{\mu}*|u|^{2^{*}_{\mu}})|u|^{2^{*}_{\mu}-2}u \ \ \mbox{in} \ \
\mathbb{R}^{N},\\ 
u>0 \ \ \mbox{in} \ \ \mathbb{R}^{N},\\
\end{array}
\right.\leqno{(P_\lambda)}
$$
where $I_\mu=\frac{1}{|x|^\mu}$ is the Riesz potential,
$0<\mu<\min\{N,4\}$ , $2^{*}_{\mu}=\frac{(2N-\mu)}{N-2}$ with $N\geq
3$, $\lambda\ge 0$ and $V_0:\mathbb{R}^N \to \mathbb{R}$ is a positive
function satisfying some technical conditions that will be mentioned
later on.    

The existence of solution for problem $(P_\lambda)$ ensures the
existence of standing waves solutions for a nonlinear Schr\"{o}dinger 
equation of the form
\begin{equation}\label{EN}
i\partial_{t} \Psi =-\Delta \Psi
+V_0(x)\Psi-(I_\mu\ast|\Psi|^{2^{*}_{\mu}})|\Psi|^{2^{*}_{\mu}-2}\Psi,\ \ \
\quad \mbox{in} \quad \mathbb{R}^N,
\end{equation}
where  $V_0$ is the external potential and $I_\mu$ is the
response function, which possesses information on the mutual interaction
between the bosons. This type of nonlocal equation appears in a lot of
physical applications, for instance in the study of  propagation of
electromagnetic 
waves in plasmas \cite{BeC} and in the theory of Bose-Einstein
condensation \cite{D}. 
We recall that a standing wave solution is a
solution of the type 
\begin{equation}
\label{an}
\Psi(x,t)=u(x)e^{{i\lambda}t},
\end{equation}
which solves
\eqref{EN} if, and only if, $u$ solves the equation 
\begin{equation}\label{SN}
-\Delta u +(\lambda+V_0(x))u
=\left(\frac{1}{|x|^{\mu}}\ast|u|^{2^{*}_{\mu}-2}\right)|u|^{2^{*}_{\mu}-2}u
\quad \mbox{in} \quad \mathbb{R}^N, 
\end{equation}
which is a Choquard-Pekar equation. 
In the ansatz (\ref{an}) the frequency $\lambda$ is related to the energy of the particle, so
$\lambda\neq 0$ is a significant case.

In \cite{DuYang}, Du and Yang have considered only the case $\lambda=V_0=0$, and
they showed that any positive solution of  (\ref{SN}) must be of the
form  
$$
\Psi_{\delta,y}(x)=C\biggl(\frac{\delta}{\delta^{2}+|x-y|^{2}}\biggl)^{\frac{N-2}{2}},
\quad x \in \mathbb{R}^N, 
$$
for some $\delta>0, y \in \mathbb{R}^N$, and $C>0$ is a constant that
depends only on $N$.  Still related to (\ref{SN}), Du, Gao and Yang
\cite{DGY} 
have  studied existence and qualitative properties of solutions of the problem
\begin{equation}\label{SN1}
-\Delta
u=\frac{1}{|x|^{\alpha}}\left(\frac{1}{|x|^{\mu}}
\ast|u|^{2^{*}_{\mu}-2}\right)|u|^{2^{*}_{\mu}-2}u 
\quad \mbox{in} \quad \mathbb{R}^N, 
\end{equation}
for some values of $\alpha$ and $\mu$. In that paper, the authors has
proved an interesting version of the Concentration-Compactness
principle due to Lions \cite{Lions} that can be used for Choquard
equations with critical growth, for more details see \cite[Lemma
2.5]{DGY}.  

In \cite{GSYZ}, Gao, da Silva, Yang and Zhou showed the
existence of solution for (\ref{SN}) by supposing $\lambda=0$ and the following
conditions on potential $V_0$: \\ 
\noindent $(I) \,\, V_0 \in C(\mathbb{R}^N,\mathbb{R}), \quad V_0(x) \geq
\nu>0$ in a neighborhood of \quad 0. \\ 
\noindent $(II)$ There are $p_1<N/2$, $p_2>N/2$ and for $N=3$, $p_2 <
3$, such that $V_0 \in L^{p}(\mathbb{R}^N), \quad \forall p \in
[p_1,p_2]$ \\ 
\noindent $(III)$
$|V_0|_{N/2}<C(N,\mu)^{\frac{N-2}{2N-\mu}}S_{H,L}(2^{\frac{N=2-\mu}{2N-\mu}}-1)$. \\ 
\noindent where $C(N,\mu)$ and $S_{H,L}$ are as in
(\ref{special}) and (\ref{melhorconstante}), respectively.

A first result fronting problems like (\ref{SN}) is due to Benci and
Cerami in the seminal paper  \cite{BC}, where the authors
studied the existence of solution for the following class of local
critical problem  
\begin{equation}\label{S.S} 
-\Delta u +V_0(x)u  = |u|^{2^{*}-2}u  \quad \mbox{in} \quad \mathbb{R}^N,
\end{equation} 
with $N \geq 3$, $2^*=\frac{2N}{N-2}$  and the function
$V_0:\mathbb{R}^N \to \mathbb{R}$ satisfies the conditions below
\begin{enumerate}
	
	\item[$(i)$]   $V_0\ge 0$ and is strictly positive in an open set.

	\item[$(ii)$]  $V_0 \in L^{q}(\mathbb{R}^{N})$ for all $q \in
          [p_{1},p_{2}]$ with $1<p_{1}<\frac{N}{2}<p_{2}$, with
          $p_2<3$ if $N=3$. 
	
	\item[$(iii)$] $|V_0|_{L^{N/2}(\mathbb{R}^{N})}< S(2^{2/N}-1),$  \\
	where $S$ denotes the best constant of the immersion $D^{1,
          2}(\mathbb{R}^{N})\hookrightarrow
        L^{2^{*}}(\mathbb{R}^{N})$, that is,  
	$$
	S:=\displaystyle\inf_{u \in D^{1,2}(\mathbb{R}^{N}),u\neq 0}
	\frac{\int_{\mathbb{R}^N}|\nabla u|^2\,dx}{|u|^{2}_{2^{*}}}.
	$$ 
\end{enumerate}
By using variational methods, the authors were able to prove the
existence of a positive solution $u \in D^{1,2}(\mathbb{R}^N)$ with 
$$
f(u) \in (S\,,\,2^{2/N}S), 
$$ 
where $f:D^{1,2}(\mathbb{R}^N)\rightarrow\mathbb{R}$ is the functional given by
$$
f(u)=\ds\int_{\mathbb{R}^N}\ds\left(|\nabla u|^{2}+V_0(x)|u|^{2}\ds\right)dx.
$$
The main difficulty to prove the existence of solution comes from the
fact that the nonlinearity has a critical growth. To overcome this
difficult, the authors used Variational methods, Deformation lemma,
and the well known  Concentration-Compactness principle due to Lions
\cite{Lions}. After the publication of \cite{BC}, some authors studied
problems related to (\ref{S.S}), see for example, \cite{Alves},
\cite{BGM1}, \cite{BNW}, \cite{CM}, \cite{Chab}, \cite{HMPY},  \cite{MW}, \cite{M},
\cite{MP}, \cite{PA}  and references therein.

In the present paper, we intend to study the existence and
multiplicity of positive solutions for $(P_{\lambda})$ for a new class of
potential $V_0$, and moreover, we also consider the case $\lambda>0$
that is a novelty for this class of problem. The assumptions on
potential $V_0$  are the following:  

\begin{enumerate}

\item[($V_{1}$)]  $V_0(x)\geq 0, \quad \forall x \in \mathbb{R}^N $,

\item[($V_{2}$)] $V_0 \in L^{N/2}(\mathbb{R}^{N})$,

\item[($V_{3}$)] $0<|V_0|_{L^{N/2}(\mathbb{R}^{N})}<
 \left(2^{\frac{4-\mu}{2N-\mu}}-1\right)\, S$.
\end{enumerate}

Our main results have the following statements: 

\begin{teorema} \label{T1} 
Let $\lambda=0$ and assume that   $(V_1)-(V_3)$ hold.
Then problem $(P_0)$ has at least a positive solution.
\end{teorema}

\begin{teorema} \label{T2} 
Let $\lambda>0$ and assume that   $(V_1)$ and $(V_2)$  hold. 
Then there exists $\bar\lambda>0$,  such that if $\lambda\in(0,\bar\lambda)$
then problem  $(P_\lambda)$ has at least a positive solutions $v_l$.
If $(V_3)$ also holds, then there exists $\lambda_0=\lambda_0(V_0)>0$
such that if  
$\lambda\in(0,\lambda_0)$ then problem  $(P_\lambda)$ has at least two
distinct positive solutions, $v_l$ and $v_h$. 
\end{teorema}

\begin{obs}
{\em
a) The solutions we find are bound state solutions, indeed
assumptions $(V_1)$ and $(V_2)$ imply that problem $(P_{\lambda})$ has no 
ground state solutions for every $\lambda\ge 0$ (see Proposition
\ref{Proposition3.1}).

b) The solution $v_h$ is an high energy one, while  $v_l$ is a low energy one.
Namely, as $\lambda\to 0$, the solution $v_h$ converges to the
solution provided by Theorem \ref{T1} while  $v_l$ flattens and
disappears. 
}\end{obs}

The existence of a solution as in Theorem \ref{T1} has been suggested
in \cite[Problem 3]{MS17} and it is proved in \cite[Theorem 1.4]{GSYZ}
with different conditions on potential $V_0$ as mentioned
above. Hence,  the Theorem \ref{T1} complements the study made
\cite{GSYZ}, in the sense that we find a positive solution for a new class of
potential $V_0$ in the case $\lambda=0$.   
Related to the Theorem \ref{T2}, we would like to point out that it is
an important contribution of this paper, because it establishes the
existence of at least two solutions for $\lambda>0$ and small enough,
which is new for this class of problem.  In order to prove the
theorems above, we will approach the problem by variational methods. 
Since the problem lacks of compactness, both for the critical exponent
and for the unboundedness of the domain, 
an analysis of the behaviour of the Palais-Smale sequences needs.
We perform it by almost classical techniques developed in the local
case, see \cite{BC}, that in our nonlocal framework requires some
careful adjustment (see Theorem \ref{Teorema de Compacidade Global}). 
In particular, the case $\lambda>0$ presents some more difficulties
because in such a case the natural space to work in is $H^1(\R^N)$
while typically the splitting theorem works in $D^{1,2}(\R^N)$. 
To front these difficulties we follow some ideas from \cite{CM} and
\cite{LM} together with a nonexistence result contained in \cite{MS}. 
Let us remark that in the critical case no extra regularity
assumptions need for the nonexistence result, see Theorem
\ref{nonET}.

Before concluding this introduction, we would like to mention that
there is a rich literature associated with Choquard-Pekar equation 
of the type 
\begin{equation}\label{SN0}
-\Delta u +V(x)u  =K(x)\left(\frac{1}{|x|^{\mu}}\ast H(u)\right)h(u)
\quad \mbox{in} \quad \mathbb{R}^N, 
\end{equation}
where $H(t)=\int_{0}^{t}h(s)ds$, with $V,K:\mathbb{R}^N \to \mathbb{R}$
and $h:\mathbb{R} \to \mathbb{R}$ being continuous functions verifying
some technical conditions.
The reader can find some interesting results in  
\cite{AGY}, \cite{AlvesNobregaYang}, \cite{AYang}, \cite{AlvesJianfu},
\cite{CCS}, \cite{L}, \cite{Lions}, \cite{ML}, \cite{MS4}, \cite{MS},
\cite{MS2}, \cite{MS3}, \cite{S} and references theirein.

The paper is organized as follows: In Section 2, we prove some results
involving the limit problem. In Section 3, we prove a splitting
theorem and show some compactness results involving the energy
functional associated with $(P_\lambda)$. In Section 4, we make the
proof of some technical lemmas that will be used in Section 5 in the
proofs of  Theorems  \ref{T1} and \ref{T2}. 

\section{Variational framework}

In this section, we will show some important results involving the limit problem that are crucial in our approach. To begin with, we recall that to apply variational methods,  we must have 
\begin{equation} \label{Z1}
 \left|\int_{\mathbb{R}^{N}}\big(I_{\mu}* |u|^{2{^{*}_\mu}}\big)  |u|^{2{^{*}_\mu}}\,dx\right|
< +\infty, \quad \forall u \in D^{1,2}(\mathbb{R}^{N}).
\end{equation}
This fact is an immediate consequence of the Hardy-Littlewood-Sobolev
inequality, which will be frequently used in this paper. 

\begin{proposicao}[\cite{LL}]
\label{Hardy} 
$\,\,[Hardy-Littlewood-Sobolev \ inequality]$:
\\
	Let $s, r>1$ and $0<\mu<N$ with $1/s+\mu/N+1/r=2$. 
If $g \in
	L^s(\mathbb{R}^N)$ and $h\in L^r(\mathbb{R}^N)$, then there exists a
	sharp constant $C(s,N,\mu,r)$, independent of $g,h$, such that
	$$
	\int_{\mathbb{R}^{N}}\int_{\mathbb{R}^{N}}\frac{g(x)h(y)}{|x-y|^{\mu}}\leq
	C(s,N,\mu,r) |g|_{s}|h|_{r}.
	$$
\end{proposicao}

As a direct consequence of this inequality, we have 
\begin{eqnarray}\label{special}
\left(\int_{\mathbb{R}^{N}}\int_{\mathbb{R}^{N}}
\frac{|u(x)|^{2^{*}_{\mu}}|u(y)|^{2^{*}_{\mu}}}{|x-y|^{\mu}} 
  dx dy\right)^{\frac{1}{2*_\mu}}\leq 
C(N,\mu)^{\frac{1}{2^*_\mu}} |u|^{2}_{2^{*}}, \quad \forall u \in
D^{1,2}(\mathbb{R}^{N}), 
\end{eqnarray}
for a suitable positive constant $C(N,\mu)$.

In the sequel, if $\lambda>0$ we will work with
$H=H^{1}(\mathbb{R}^{N})$ endowed with the norm 
$$
\|u\|=\left(\displaystyle\int_{\mathbb{R}^{N}}|\nabla u|^{2} dx+
\lambda\displaystyle\int_{\mathbb{R}^{N}}|u|^{2} dx\right)^{\frac{1}{2}}.   
$$
When $\lambda =0$, we will consider $H=D^{1,2}(\mathbb{R}^{N})$ endowed with the usual norm, that is, 
$$
\|u\|=\left(\displaystyle\int_{\mathbb{R}^{N}}|\nabla u|^{2} dx\right)^{\frac{1}{2}}. 
$$
Sometimes we use the explicit notations $\|\cdot\|_{D^{1,2}}$ and $\|\cdot\|_{H^1}$.

We say that $u:\mathbb{R}^N\rightarrow\mathbb{R}$ is a weak solution
of $(P_\lambda)$ if $u\in H$ is a positive function such that for all
$\varphi\in H$ we get 
$$
\displaystyle\int_{\mathbb{R}^N}\nabla u \nabla\varphi
dx+\displaystyle\int_{\mathbb{R}^N}V_{\lambda}(x)u\varphi
dx=\displaystyle\int_{\mathbb{R}^N}(I_{\mu}*|u|^{2^{*}_{\mu}})|u|^{2^{*}_{\mu}-2}u\varphi
dx, 
$$
where $V_{\lambda}:=\lambda+V_{0}$. 
It is a standard task to check that the weak
solutions of $(P_\lambda)$ are critical points of the energy
functional  
$E_{\lambda}:H\to \mathbb{R}$ associated to problem $(P_{\lambda})$ 
given by
$$
E_{\lambda}(u)=\frac{1}{2}\|u\|^{2} +
\frac{1}{2}\displaystyle\int_{\mathbb{R}^{N}}V_{0}(x)u^{2} dx -
\frac{1}{2\cdot
  2^{*}_{\mu}}\displaystyle\int_{\mathbb{R}^N}
(I_{\mu}*|u|^{2^{*}_{\mu}})|u|^{2^{*}_{\mu}}dx,  
$$
and that  $E_{\lambda}$ belongs to $C^{1}(H,\mathbb{R})$.

In what follows, let us denote by $\mathcal{M}$ the manifold 
$$
\mathcal{M}=\ds\bigg\{u\in H\ :\
\ds\int_{\mathbb{R}^N}(I_{\mu}*|u|^{2^{*}_{\mu}})|u|^{2^{*}_{\mu}}dx=1\ds\bigg\} 
$$
and by  $J_{\lambda}$ the functional
$$
J_{\lambda}(u)=\|u\|^{2} + \displaystyle\int_{\mathbb{R}^{N}}V_{0}(x)
u^{2} dx,\quad\qquad u\in H.
$$
For the case $\lambda=0$, we designate such functionals by $E_{0}$ and
$J_{0}$, respectively. 

Next, let us recall some information involving the problem  
$$
 \left\{
\begin{array}{rcl}
-\Delta u=(I_{\mu}*|u|^{2^{*}_{\mu}}) |u|^{2^{*}_{\mu}-2}u \ \ \mbox{in} \ \ \mathbb{R}^{N},\\
u \in D^{1,2}(\mathbb{R}^N).
\end{array}
\right.\leqno{(P_{\infty})}
$$
The functions
\begin{equation}  \label{tildeU}
\tilde{U}_{\sigma,z}(x)=
\frac{[N(N-2)\sigma]^{\frac{N-2}{4}}}{(\sigma+|x-z|^{2})^{\frac{N-2}{2}}}, 
\quad \sigma>0,  \quad z\in \mathbb{R}^{N}, 
\end{equation}
are the minimizers for $S$ in $D^{1,2}(\R^N)$ that verify  
$$
 \left\{
\begin{array}{rcl}
-\Delta u=|u|^{2^{*}-2}u \ \ \mbox{in} \ \ \mathbb{R}^{N},\\
u \in D^{1,2}(\mathbb{R}^N)
\end{array}
\right.
$$
and 
\begin{equation}  \label{U}
U_{\sigma,z}(x)=C(N,\mu)^{\frac{2-N}{2(N-\mu+2)}}
S^{\frac{(N-\mu)(2-N)}{4(N-\mu+2)}}\tilde{U}_{\sigma,z}(x),
\quad \forall x \in \mathbb{R}^N 
\end{equation}
are minimizer  for $S_{H,L}$ given by 
\begin{eqnarray}\label{melhorconstante}S_{H,L}:=\displaystyle\inf_{u
    \in D^{1,2}(\mathbb{R}^{N}),u\neq 0} 
\frac{\|u\|^{2}}{\biggl(\displaystyle\int_{\mathbb{R}^{N}}
(I_{\mu}*|u|^{2^{*}_{\mu}})|u|^{2^{*}_{\mu}}
  dx\biggl)^{1/2^{*}_{\mu}}},  
\end{eqnarray}
that satisfy $(P_{\infty})$ (see \cite{Talenti,L83} and \cite[Theorem
4.3]{LL}), with  
\begin{eqnarray}\label{vouusar}
\displaystyle\int_{\mathbb{R}^{N}}|\nabla U_{\sigma,z}|^{2} dx = 
\int_{\mathbb{R}^N}(I_{\mu}*|U_{\sigma,z}|^{2^{*}_{\mu}})|U_{\sigma,z}|^{2^{*}_{\mu}}dx=
S_{H,L}^{\frac{2N-\mu}{N+2-\mu}}.
\end{eqnarray}
Moreover,  it is possible to show that
$$
C(N,\mu)^{\frac{1}{2^*_\mu}}S_{H,L}=S \quad (\mbox{See} \,\,
\cite[Lemma \, 1.2]{Gao}) 
$$
and
\beq
\label{minimodeI_0}
E_{\infty}(U_{\sigma,z})=\frac{(N+2-\mu)}{4N-2\mu}S_{H,L}^{(2N-\mu)/(N+2-\mu)}, 
\eeq
where $E_{\infty}$ is the energy functional associated to $(P_{\infty})$ given by
$$
E_{\infty}(u)=\frac 12 \|u\|^{2} - \frac{1}{2\cdot 2^*_\mu}
\int_{\mathbb{R}^{N}}\big(I_{\mu}* |u|^{2{^{*}_\mu}}\big)
|u|^{2{^{*}_\mu}}\,dx. 
$$

According to the limit problem when $\lambda>0$, we have the result
below that is a key point in our approach. 
\begin{teorema}
\label{nonET}
If $\lambda>0$ and $u\in H^1(\R^N)$  is a weak solution of
$$
 -\Delta u+\lambda u = (I_{\mu}*|u|^{2^{*}_{\mu}})|u|^{2^{*}_{\mu}-2}u,
\leqno{(P_{\lambda,\infty})}
$$
then $u\equiv 0$.
\end{teorema}
\noindent {\bf Proof.}
In what follows, we are proving the assumption of \cite[Theorem 2]{MS} holds. 
We have only to verify that in our critical case all the regularity
required in the proof is fulfilled. Since here $p=2^*_\mu$, so
$\frac{1}{p}=\frac{N-2}{N+\alpha}$ where $\alpha=N-\mu$, it is
sufficient to show that if $u\in H^1(\R^N)$  is a weak solution of
$(P_{\infty,\lambda})$ then $u\in H^{2,2}_{\loc}(\R^N)$. 
The function $u$ solves 
$$
-\Delta u=(-\lambda +Q(x))u \quad u\in H^1(\mathbb{R}^N), 
$$
where
$$
Q(x)=
\Big( I_\mu\ast |u|^{2^*_\mu}\Big)|u|^{2_\mu^*-2}.
$$
By the Hardy-Littlewood-Sobolev inequality, we know that $I_\mu\ast
|u|^{2^*_\mu} \in L^{\frac{2N}{\mu}}(\mathbb{R}^N).$ Since
$|u|^{2^*_\mu-2} \in L^{\frac{2N}{4-\mu}}(\mathbb{R}^N) $, we claim
that $Q \in L^{\frac{N}{2}}(\mathbb{R}^N)$. 
Indeed, as
$$
|Q(x)|^{\frac{N}{2}}=(I_\mu\ast |u|^{2^*_\mu})^{\frac{N}{2}}(|u|^{2^*_\mu-2})^{\frac{N}{2}}
$$
and
$$
(I_\mu\ast |u|^{2^*_\mu})^{\frac{N}{2}} \in
L^{\frac{4}{\mu}}(\mathbb{R}^N), \quad (|u|^{2^*_\mu-2})^{\frac{N}{2}}
\in L^{\frac{4}{4-\mu}}(\mathbb{R}^N), 
$$
we derive from the H\"{o}lder inequality that $|Q|^{\frac{N}{2}}$
belongs to $L^{1}(\mathbb{R}^N)$, that is, $Q \in
L^{\frac{N}{2}}(\mathbb{R}^N)$. 
Now, arguing as in Struwe \cite[Lemma B.3]{StruweBook},
$$
u \in L_{loc}^{q}(\mathbb{R}^N), \quad \forall q \in[1,+\infty).
$$
Now, we claim that $I_\mu\ast |u|^{2^*_\mu} \in
L^{\infty}(\mathbb{R}^N)$. Indeed, note that,  
$$\aligned
|(I_\mu\ast |u|^{2^*_\mu})(x)|&=\int_{\R^N}\frac{|u|^{2^*_\mu}}{|x-y|^\mu}dy\\
&\leq \int_{|y-x|\leq1}\frac{|u|^{2^*_\mu}}{|x-y|^\mu}dy+
\int_{|y-x|\geq1}\frac{|u|^{2^*_\mu}}{|x-y|^\mu}dy.\\ 
\endaligned
$$
Using the fact that
$$
\frac{1}{|z|^{\mu}}\in L^{\frac{2N}{\mu}}(B^{c}_1(0)),  \quad
|u|^{2^*_\mu} \in L^{\frac{2N}{2N-\mu}}(B^{c}_1(x)),\quad
\frac{1}{|z|^{\mu}}\in L^{t}(B_1(0)) \quad \mbox{for} \quad t \in [1,
{N}/{\mu})  
$$
and for every $x\in\R^N$
$$
|u|^{2^*_\mu} \in L^{t'}(B_1(x)) \quad \mbox{for }\ t'\ \mbox{ s.t. }
\quad\frac{1}{t} + 
\frac{1}{t'} =1, 
$$
the H\"{o}lder inequality ensures that for every compact set
$K\subset\R^N$ there is $C=C(K)>0$ such that
$$
|(I_\mu\ast |u|^{2^*_\mu})(x)|\leq C, \quad \forall x \in K,
$$
which implies $I_\mu\ast |u|^{2^*_\mu} \in
L^{\infty}_{\loc}(\mathbb{R}^N)$. 
Recalling that $|u|^{2^*_\mu-2} \in  L_{loc}^{q}(\mathbb{R}^N)$ for
every $q \in [1,+\infty)$, it follows that  
$$
Q \in  L_{loc}^{q}(\mathbb{R}^N) \quad \forall q \in [1,+\infty).
$$
Thereby, by Calder\'on-Zygmund inequality, see \cite{GT}, $u \in
W^{2,q}_{loc}(\mathbb{R}^N)$ for all $q \in (1,+\infty)$. This proves
the theorem. 
{\fim}

Before concluding this section, we will prove an important estimate
involving the nodal solutions of the limit problem, which will be used
later on.     
\begin{lema} \label{Nodal}
If $u \in D^{1,2}(\mathbb{R}^N)$ is a nodal solution of $(P_{\lambda})$, then
$$
E_{\lambda}(u) \geq
2^{\frac{4-\mu}{N+2-\mu}}\frac{(N+2-\mu)}{4N-2\mu}S_{H,L}^{(2N-\mu)/(N+2-\mu)}. 
$$	
\end{lema} 
\noindent {\bf Proof.} 
Arguing as in the proof of \cite[Proposition 3.2]{GV1}, for all
$t^{+},t^->0$ we see that 
$$
E_{\lambda}(t^+u^+)^{\frac{N+2-\mu}{4-\mu}} +
E_{\lambda}(t^-u^{-})^{\frac{N+2-\mu}{4-\mu}} \leq
E_{\lambda}(u)^{\frac{N+2-\mu}{4-\mu}}.  
$$
Fixing $t^+,t^- >0$ such that
$E'_{\lambda}(t^{\pm}u^{\pm})(t^{\pm}u^{\pm})=0$, it follows that 
$$
E_{\lambda}(t^{\pm}u^{\pm}) \geq \frac{(N+2-\mu)}{4N-2\mu}S_{H,L}^{(2N-\mu)/(N+2-\mu)}.
$$
The last two inequalities combine to give 
$$
E_{\lambda}(u) \geq
2^{\frac{4-\mu}{N+2-\mu}}\frac{(N+2-\mu)}{4N-2\mu}S_{H,L}^{(2N-\mu)/(N+2-\mu)}, 
$$
finishing the proof. 

{\fim}

\begin{lema}\label{SHLS}
If $(u_n) \subset D^{1,2}(\mathbb{R}^N)$ is such that 
$$
\frac{\displaystyle \int_{\mathbb{R}^{N}}|\nabla u_n|^{2}
  dx}{\biggl(\ds\int_{\mathbb{R}^N}(I_{\mu}*|u_n|^{2^{*}_{\mu}})|u_n|^{2^{*}_{\mu}})dx
\biggl)^{1/2^{*}_{\mu}}}
\to S_{H,L}, 
$$
then there are $\sigma_n>0$ and $y_n \in \mathbb{R}^N$ such
that
$$
\frac{u_n}{|u_n|_{2^*}}=\frac{{U}_{\sigma_n,y_n}}{|{U}_{\sigma_n,y_n}|_{2^*}}+o_n(1)
\quad \mbox{in} \quad D^{1,2}(\mathbb{R}^N), 
$$
where ${U}_{\sigma,y}$ is given in (\ref{U}).
\end{lema}
\noindent {\bf Proof.}  Note that by (\ref{special}), 
$$
S_{H,L}+o_n(1) =\frac{\displaystyle\int_{\mathbb{R}^{N}}|\nabla
  u_n|^{2}
  dx}{\biggl(\ds\int_{\mathbb{R}^N}
(I_{\mu}*|u_n|^{2^{*}_{\mu}})|u_n|^{2^{*}_{\mu}}dx\biggl)^{1/2^{*}_{\mu}}}\geq 
\frac{\displaystyle \int_{\mathbb{R}^{N}}|\nabla u_n|^{2}
  dx}{C(N,\mu)^{\frac{1}{2^*_\mu}} |u_n|^{2}_{2^{*}}}\geq
\frac{S}{C(N,\mu)^{\frac{1}{2^*_\mu}}}. 
$$
Since $S_{H,L}C(N,\mu)^{\frac{1}{2^*_\mu}}=S$, then 
$$
\frac{\displaystyle \int_{\mathbb{R}^{N}}|\nabla u_n|^{2}
  dx}{\biggl(\ds\int_{\mathbb{R}^N}|u_n|^{2^{*}}dx\biggl)^{2/2^{*}}}
\to S. 
$$
By Lions \cite{Lions}, Aubin \cite{Aubin} and Talenti \cite{Talenti},
there are ${\sigma}_n>0$ and ${y}_n \in \mathbb{R}^N$ such
that  
$$
\frac{u_n}{|u_n|_{2^*}}=
\frac{\tilde{U}_{{\sigma}_n,{y}_n}}
{|\tilde{U}_{{\sigma}_n,{y}_n}|_{2^*}}+o_n(1) 
\quad \mbox{in} \quad D^{1,2}(\mathbb{R}^N), 
$$
where $\tilde{U}_{\sigma,y}$ is given in (\ref{tildeU}). This together with (\ref{U}) proves the desired result. {\fim}

Next corollary states a first property of Palais-Smale sequences. 
For short, we shall refer to those sequences as (PS$)_c$
sequences, where $c$ is the energy level, or as (PS) sequences.

\begin{corolario} \label{compacidadeB} Let $(u_n)$ be a   $(PS)$
  sequence  for the  functional 
	$E_{\infty}$ with 
	$$
	c=\frac{(N+2-\mu)}{4N-2\mu}S_{H,L}^{(2N-\mu)/(N+2-\mu)}.
	$$
	Then, there are sequences $(\sigma_n)\subset\mathbb{R}^+$ and
	$(y_n)\subset\mathbb{R}^{N}$  such that, up to a subsequence of
	$(u_n)$, we have 
	\begin{eqnarray}
	\frac{u_n}{|u_n|_{2^*}}=\frac{{U}_{\sigma_n,y_n}}{|{U}_{\sigma_n,y_n}|_{2^*}}+o_n(1),
	\ \ \mbox{in} \ \ D^{1,2}(\mathbb{R}^N).
	\end{eqnarray}
	
\end{corolario}
\noindent {\bf Proof.}  
Since $(u_n)$ is a $(PS)_c$ sequence for $E_\infty$ it is bounded and
so it is possible to prove that 
$$
\frac{\displaystyle \int_{\mathbb{R}^{N}}|\nabla u_n|^{2}
	dx}{\biggl(\ds\int_{\mathbb{R}^N}
	(I_{\mu}*|u_n|^{2^{*}_{\mu}})|u_n|^{2^{*}_{\mu}})dx\biggl)^{1/2^{*}_{\mu}}}
\to S_{H,L}. 
$$
Now, we apply Lemma \ref{SHLS} and, taking into account that  the
weak limit of $u_n$ (up to a rescaling) solves  $(P_\infty)$, we get
the desired result.  
{\fim}

\begin{lema}
\label{PSbdd}
Let $(u_n)$ be a (PS) sequence for $E_\lambda$. 
Then $(u_n)$ is bounded and there exists a weak solution $u_0$ of
$(P_\lambda)$ such that $u_n \rightharpoonup u_0$ in $H$, up to a
subsequence. The same statement holds for (PS) sequences for the
functional $E_\infty$, with $u_0$ in $D^{1,2}(\R^N)$ a solution of
$(P_\infty)$. 
\end{lema}
\noindent {\bf Proof.} 
Since $V_0\ge 0$ on $\R^N$, then
\begin{eqnarray*}
c+o_n(1) \|u_n\|& = & E_\lambda(u_n)-\frac{1}{2\cdot 2^*_\mu}\,
E'_\lambda(u_n)u_n\\
 &\ge &
\frac{N+2-\mu}{4N-2\mu}\, \|u_n\|^2,
\end{eqnarray*}
so we get $(\|u_n\|)$ bounded and the existence of the weak limit
$u_0$.
Moreover, clearly $u_0$ is a weak solution of $(P_\lambda)$ because
$(u_n)$ is a (PS) sequence for $E_\lambda$. The same argument works for $E_\infty$.

{\fim}

\begin{lema} \label{CH} Let $\Omega \subset \mathbb{R}^N$ a smooth
  bounded domain. For each $\tau>0$, there is $M_\tau>0$ such that  
	$$
	\left[\frac{S}{2^{\frac{2}{N}}}-\tau
        \right]^{2^*_\mu}\int_{\Omega}(I_{\mu}*|v|^{2^*_{\mu}})|v|^{2^*_\mu}\,dx
        \leq C(N,\mu) \left(|\nabla
          v|^{2}_{L^{2}(\Omega)}+M_\tau|v|_{L^{2}(\Omega)}^{2}\right)^{2^*_\mu},
        \quad \forall v \in H^{1}(\Omega), 
	$$
	where $C(N,\mu)$ is given in (\ref{special}).
\end{lema}
\noindent {\bf Proof.} From (\ref{special}), 
$$
\left[\frac{S}{2^{\frac{2}{N}}}-\tau
\right]^{2^*_\mu}\int_{\Omega}(I_{\mu}*|v|^{2^*_{\mu}})|v|^{2^*_\mu}\,dx
\leq C(N,\mu) \left(\left[\frac{S}{2^{\frac{2}{N}}}-\tau
  \right]|v|^2_{L^{2^*}(\Omega)}\right)^{2^*_\mu}. 
$$
Now, we apply the Cherrier's inequality \cite{C} to get the desired result.
{\fim}

\section{A global compactness theorem}

In this section we study the Palais-Smale sequences of our functionals
First,  now prove a technical lemma for $E_\infty$ that will be useful later on. 

\begin{lema}(Main lemma) \label{Lematecnico}
Let $(u_n)$ be a   $(PS)_{c}$ sequence  for the  functional
$E_{\infty}$ with $u_n\rightharpoonup 0$ and $u_n\nrightarrow
0$. Then, there are sequences $(\sigma_n)\subset\mathbb{R}^+$ and
$(y_n)\subset\mathbb{R}^{N}$  such that, up to a subsequence of
$(u_n)$, we have 
\begin{eqnarray}
	v_n(x): =\sigma^{(N-2)/2}_{n} u_n(\sigma_{n}x+y_n) \rightharpoonup v_0,
	\ \ \mbox{in} \ \ D^{1,2}(\mathbb{R}^N) 
\end{eqnarray}
where $v_0$ is a non trivial solution of problem $(P_{\infty})$ and
\begin{eqnarray}\label{conv2}
	E_{\infty}(u_n)= E_{\infty}(v_0)+E_{\infty}(v_n-v_0)+o_{n}(1)
\end{eqnarray}
\begin{eqnarray}\label{conv3}
	\|u_n\|^{2}=\|v\|^{2}+ \|v_n -v\|^{2}+o_{n}(1).
\end{eqnarray}
\end{lema}

\noindent {\bf Proof.}   
Let $(u_n)\subset D^{1,2}(\mathbb{R}^{N})$ be a $(PS)_{c}$ sequence
for the functional $E_{\infty}$, i.e, 
\begin{eqnarray}\label{LP1}
E_{\infty}(u_n)\rightarrow c\,\,\ {\rm and}\,\,\ E'_{\infty}(u_n)\rightarrow 0.
\end{eqnarray}
The sequence $(u_n)$ is bounded in $D^{1,2}(\mathbb{R}^{N})$, by Lemma \ref{PSbdd}. 
Then, since $u_n\rightharpoonup 0$ and $u_n\nrightarrow 0$, the (PS) condition and (\ref{melhorconstante}) imply that
$$
c\geq\frac{(N+2-\mu)}{4N-2\mu}S_{H,L}^{(2N-\mu)/(N+2-\mu)}.
$$
Note that
\begin{eqnarray*}
c +o_{n}(1)=E_{\infty}(u_n)-\displaystyle\frac{1}{2\cdot 2^{*}_{\mu}}E'_{\infty}(u_n)u_n =\frac{(N+2-\mu)}{4N-2\mu}\displaystyle\int_{\mathbb{R}^{N}}|\nabla u_n|^{2}dx
\end{eqnarray*}
which leads to
\begin{eqnarray}\label{LP2}
\displaystyle\lim_{n\rightarrow \infty}\displaystyle\int_{\mathbb{R}^{N}}|\nabla u_n|^{2}dx\geq S_{H,L}^{(2N-\mu)/(N+2-\mu)}.
\end{eqnarray}
Let $m$ be a number such that $B_{2}(0)$  is covered by $m$ balls of radius $1$, $k \in \mathbb{N}^*$ to be fixed later, $(\sigma_n)\subset\mathbb{R}^+$, $(y_n)\subset\mathbb{R}^{N}$ such that 
$$
\displaystyle\sup_{y\in\mathbb{R}^{N}}\displaystyle\int_{B_{\sigma_{n}}(y)}|\nabla u_n|^{2}dx=\displaystyle\int_{B_{\sigma_{n}}(y_n)}|\nabla u_n|^{2}dx=\displaystyle\frac{S_{H,L}^{(2N-\mu)/(N+2-\mu)}}{km}
$$
and 
$$
v_{n}(x):=\sigma^{(N-2)/2}_{n}u_n\displaystyle\bigg(\sigma_{n}x+y_{n}\displaystyle\bigg).
$$
A simple computation gives 
$$
\|v_n\|=\|u_n\| \quad \mbox{and} \quad \ds\int_{\mathbb{R}^{N}}(I_{\mu}*|v_n|^{2^{*}_{\mu}})|v_n|^{2^{*}_{\mu}}dx=\ds\int_{\mathbb{R}^{N}}(I_{\mu}*|u_n|^{2^{*}_{\mu}})|u_n|^{2^{*}_{\mu}}dx,
$$
from where it follows that 
\begin{equation} \label{IG1}
E_\infty(v_n)=E_\infty(u_n), \quad \forall n \in \mathbb{N}.
\end{equation}

Hereafter, we fix $k$ such way that  
\begin{equation} \label{m}
\displaystyle\frac{S_{H,L}^{(2N-\mu)/(N+2-\mu)}}{k} <
\left(\frac{S}{2^{\frac{2}{N}}}\right)^{\frac{2N-\mu}{N+2-\mu}}
\left(\frac{1}{C(N,\mu)}\right)^{\frac{N-2}{N+2-\mu}}. 
\end{equation}
Using a change of variable, we can prove that 
\begin{equation} \label{NOVA}
\displaystyle\int_{B_{1}(0)}|\nabla
v_n|^{2}dx=\displaystyle\frac{S_{H,L}^{(2N-\mu)/(N+2-\mu)}}{km}=
\displaystyle\sup_{y\in\mathbb{R}^{N}}\int_{B_{1}(y)}|\nabla
v_n|^{2}dx. 
\end{equation}
Now, for each  $\Phi\in D^{1,2}(\mathbb{R}^{N})$, we define the function 
$$
\tilde{\Phi}_{n}(x)=\frac{1}{\sigma^{(N-2)/2}_{n}}\Phi(\frac{1}{\sigma_n}(x-y_n))
$$
that satisfies
\begin{eqnarray}\label{PP1}
\displaystyle\int_{\mathbb{R}^{N}}\nabla u_n \nabla \tilde{\Phi}_{n}dx=\displaystyle\int_{\mathbb{R}^{N}}\nabla v_n\nabla \Phi dx
\end{eqnarray}
and
\begin{eqnarray}\label{PP2}
\displaystyle\int_{\mathbb{R}^{N}}(I_{\mu} * |u_n|^{2^{*}_{\mu}})|u_n|^{2^{*}_{\mu}-2}u_n\tilde{\Phi}_{n}dx=\displaystyle\int_{\mathbb{R}^{N}}(I_{\mu} * |v_n|^{2^{*}_{\mu}})|v_n|^{2^{*}_{\mu}-2}v_n\Phi dx.
\end{eqnarray}
These limits ensure that 
\begin{equation} \label{IG2}
E'_\infty(v_n)\to 0.
\end{equation} 
From (\ref{IG1})-(\ref{IG2}), $(v_n)$ is a $(PS)_c$ sequence for $E_\infty$. Therefore, there exists $v_0\in D^{1,2}(\mathbb{R}^{N})$ such that, up to a subsequence,  $v_n\rightharpoonup v_0$ in $D^{1,2}(\mathbb{R}^{N})$ and $E'_{\infty}(v_0)=0$.

As a consequence of the well known Lions' Lemma \cite{Lions}, we can assume that
\begin{eqnarray}
\label{LP6}
\displaystyle\int_{\mathbb{R}^{N}}|v_n|^{2^{*}}\phi dx\rightarrow
\displaystyle\int_{\mathbb{R}^{N}}|v_0|^{2^{*}}\phi
dx+\displaystyle\sum_{j\in J}\phi(x_j)\nu_{j},\,\,\ \forall \phi\in
C^{\infty}_{0}(\mathbb{R}^{N}) 
\end{eqnarray}
and
$$
|\nabla v_n|^{2}\rightharpoonup \mu \geq |\nabla
v_0|^{2}+\displaystyle\sum_{j\in J}\delta_{x_j}\mu_{j}, 
$$
for some  $\{x_j\}_{j\in J}\subset\mathbb{R}^{N}$, $\{\nu_{j}\}_{j\in
  J}$, $\{\mu_{j}\}_{j\in J}\subset\mathbb{R}^{+}$ with
$S\nu_{j}^{2/2^{*}}\leq \mu_{j}$, where $J$ is at most a countable
set. 

We are going to show that $J$ is finite. 
Consider $\phi \in C^{\infty}_{0}(\mathbb{R}^{N},[0,1])$ such that  
$\phi(x)=0$ for all $x
\in B^{c}_{2}(0)$ and  $\phi_{\rho}(x)=1$ for all $x \in
B_{1}(0)$. 
Now fix $x_j \in \mathbb{R}^{N}$, $j \in J$, and define
$\psi_{\rho}(x)=\phi(\frac{x-x_j}{\rho})$, for each $\rho>0$. 
Then $0\leq \psi_{\rho}(x)\leq 1$, for all $x \in \mathbb{R}^{N}$,
$\psi_\rho(x)=0$ for all $x \in B^{c}_{2\rho}(x_j)$ and  $\psi(x)=1$
for all $x \in B_{\rho}(x_j)$. We have that $(v_{n}\psi_{\rho})$ is
bounded in $D^{1,2}(\mathbb{R}^{N})$ and
$I'_{\infty}(v_n)v_{n}\psi_{\rho} = o_{n}(1)$. 
Hence,
\begin{eqnarray}\label{vaiusar} 
\displaystyle\int_{\mathbb{R}^{N}}\psi_{\rho}|\nabla v_n|^{2} dx +
\displaystyle\int_{\mathbb{R}^{N}}v_n \nabla v_n \nabla \psi_{\rho} dx
=
\displaystyle\int_{\mathbb{R}^{N}}(I_{\mu}*|v_{n}|^{2^{*}_\mu})|v_{n}|^{2^{*}_\mu}
\psi_{\rho} dx+o_{n}(1). 
\end{eqnarray}
Using Proposition \ref{Hardy}, and seeing that 
$$
\displaystyle\lim_{\rho \to 0}\biggl[\limsup_{n \to \infty}\int_{\mathbb{R}^{N}}v_n \nabla v_n \nabla \psi_{\rho} dx\biggl]=0,
$$
we find
$$
S \nu_{j}^{2/2^{*}}\leq \mu_{j} \leq C \nu_{j}^{\frac{2^{*}_{\mu}}{2^*}}.
$$
As $2^{*}_{\mu}>2$ we deduce that $\nu_{j}$ does not converge to zero 
 and since $\displaystyle{\sum_{j \in J} \nu_{j}^{2/2^{*}}\le \sum_{j \in J} \mu_{j}< +\infty}$ 
 we have that  $J$ is finite. 
 From now on, we denote by $J=\{1,2,...,m\}$ and $\Gamma\subset\mathbb{R}^{N}$ the set given by
$$
\Gamma=\{x_j\in \{x_j\}_{j\in J}\ :\  |x_j|>1\}, \,\,\ (x_j\,\,\ {\rm given \ \  by}\,\,\ (\ref{LP6})).
$$
In the sequel, we are going to show that  $v_0\neq 0$. Suppose, by contradiction that $v_0=0$. Thereby, by $(\ref{LP6})$, 
\begin{eqnarray}\label{LP7}
\displaystyle\int_{\mathbb{R}^{N}}|v_n|^{2^{*}}\tilde\varphi\, dx\rightarrow 0,\,\,\ \forall\tilde \varphi\in C_{0}(\mathbb{R}^{N}\setminus\{x_1,x_2,...,x_m\}).
\end{eqnarray}
Using again Proposition \ref{Hardy}, for all $  \varphi\in C_{0}(\mathbb{R}^{N}\setminus\{x_1,x_2,...,x_m\})$ we derive the inequality below
\begin{eqnarray}\label{sacadaClaudianor}
\displaystyle\int_{\mathbb{R}^{N}}(I_{\mu}*|v_{n}|^{2^*_\mu})|v_{n}|^{2^*_\mu} \varphi \, dx \leq C |v_n|^{2^*_\mu}_{2^*}\left(\int|v_n|^{2^*}|\varphi|^{\frac{2N}{2N-\mu}}\,dx\right)^{\frac{2N-\mu}{2N}}=o_{n}(1), 
\end{eqnarray}
which leads to 
\begin{eqnarray}\label{Gaetano}
\displaystyle\int_{\mathbb{R}^{N}}|\nabla v_n|^{2} \varphi= o_{n}(1).
\end{eqnarray}
Consequently, if $\rho\in\mathbb{R}$ is a number that satisfies  
$0<\rho<\displaystyle\min\{{\rm dist}(\Gamma,\bar{B}_{1}(0)),1)\}$, it follows that 
\begin{eqnarray}\label{Caramba}
\displaystyle\int_{B_{1+\frac{2\rho}{3}}(0)\setminus B_{1+\frac{\rho}{3}}(0)}|\nabla v_n|^{2}dx= o_{n}(1).
\end{eqnarray}
In the sequel, let us consider the sequence $(\Phi_{n})$ given by $\Phi_{n}(x)=\Phi(x)v_n(x)$, where $\Phi\in C^{\infty}_{0}(\mathbb{R}^{N},[0,1])$ satisfies $\Phi(x)=1$ if $x\in B_{1+\rho/3}(0)$ and $\Phi(x)=0$ if $x\in B^{c}_{1+2\rho/3}(0)$.
Note that,
\begin{eqnarray*}
\displaystyle\int_{B_{1+\rho}(0)\setminus B_{1+\frac{\rho}{3}}}|\nabla \Phi_{n}|^{2}dx
\leq 
C\biggl[ \displaystyle\int_{{B_{1+\frac{2\rho}{3}}}(0)\setminus B_{1+\frac{\rho}{3}}}|\nabla v_{n}|^{2}dx +  
\displaystyle\int_{B_{1+\rho}(0)\setminus B_{1+\frac{\rho}{3}}(0)} |v_{n}|^{2}dx \biggl],
\end{eqnarray*}
that is,
\begin{eqnarray}\label{Gabriel}
\displaystyle\int_{B_{1+\rho}(0)\setminus B_{1+\frac{\rho}{3}}}|\nabla \Phi_{n}|^{2}dx = o_{n}(1).
\end{eqnarray}
Since $I'_{\infty}(v_n)\Phi_{n}= o_{n}(1)$, we have
\begin{eqnarray*}
\begin{aligned}
&\displaystyle\int_{B_{1+\rho}(0)\setminus B_{1+\frac{\rho}{3}}(0)}\nabla v_n\nabla \Phi_{n}dx+\displaystyle\int_{B_{1+\frac{\rho}{3}}(0)}\nabla v_n\nabla \Phi_{n}dx\\
&- \displaystyle\int_{B_{1+\rho}(0)\setminus B_{1+\frac{\rho}{3}}(0)}(I_{\mu}*|v_n|^{2^{*}_{\mu}})|v_n|^{2^{*}_{\mu}}\Phi\, dx-\displaystyle\int_{B_{1+\frac{\rho}{3}}(0)}(I_{\mu}*|v_n|^{2^{*}_{\mu}})|v_n|^{2^{*}_{\mu}}\Phi\, dx=o_{n}(1),
\end{aligned}
\end{eqnarray*}
which implies
\begin{eqnarray}\label{LP12}
\begin{aligned}
&\displaystyle\int_{B_{1+\rho}(0)\setminus B_{1+\frac{\rho}{3}}(0)}\nabla v_n\nabla \Phi_{n}dx+\displaystyle\int_{B_{1+\frac{\rho}{3}}(0)}|\nabla v_n|^{2}dx\\
&- \displaystyle\int_{B_{1+\rho}(0)\setminus B_{1+\frac{\rho}{3}}(0)}(I_{\mu}*|v_n|^{2^{*}_{\mu}})|v_n|^{2^{*}_{\mu}}\Phi \, dx-\displaystyle\int_{B_{1+\frac{\rho}{3}}(0)}(I_{\mu}*|v_n|^{2^{*}_{\mu}})|v_n|^{2^{*}_{\mu}}\Phi\, dx=o_{n}(1).
\end{aligned}
\end{eqnarray}
Note that from H\"older's inequality and  $(\ref{Caramba})$ 
\begin{eqnarray}\label{LP13}
\displaystyle\int_{B_{1+\rho}(0)\setminus B_{1+\frac{\rho}{3}}(0)}\nabla v_n\nabla\Phi_{n}dx\rightarrow 0\,\,\ {\rm when}\,\,\ n\rightarrow\infty
\end{eqnarray}
and that (\ref{LP7}) together with Proposition \ref{Hardy} gives 
\begin{eqnarray}\label{LP14}
\displaystyle\int_{B_{1+\rho}(0)\setminus B_{1+\frac{\rho}{3}}(0)}(I_{\mu}*|v_n|^{2^{*}_{\mu}})|v_n|^{2^{*}_{\mu}}\Phi dx=o_{n}(1).
\end{eqnarray}
Thereby, from $(\ref{LP12})$, $(\ref{LP13})$ and  $(\ref{LP14})$,
\begin{eqnarray}\label{LP15}
\displaystyle\int_{B_{1+\frac{\rho}{3}}(0)}|\nabla v_{n}|^{2}dx-\displaystyle\int_{B_{1+\frac{\rho}{3}}(0)}(I_{\mu}*|v_n|^{2^{*}_{\mu}})|v_{n}|^{2^{*}_{\mu}}dx=o_{n}(1).
\end{eqnarray}
The last equality together with the boundedness of $(v_n)$ and (\ref{NOVA}) implies that for some subsequence  
$$
\lim_{n\rightarrow\infty}\int_{B_{1+\frac{\rho}{3}}(0)}|\nabla v_{n}|^{2}dx=\lim_{n\rightarrow\infty}\displaystyle\int_{B_{1+\frac{\rho}{3}}(0)}(I_{\mu}*|v_n|^{2^{*}_{\mu}})|v_{n}|^{2^{*}_{\mu}}dx=A>0.
$$
These limits combined with the Cherrier's inequality (see Lemma \ref{CH}) give 
\begin{equation} \label{A}
A \geq \left(\frac{S}{2^{\frac{2}{N}}}\right)^{\frac{2N-\mu}{N+2-\mu}}\left(\frac{1}{C(N,\mu)}\right)^{\frac{N-2}{N+2-\mu}}.
\end{equation} 
On the other hand,
since $B_{1+\frac{\rho}{3}(0)}\subset B_{2}(0)$ and $B_2(0)$ is covered by $m$ balls of radius 1, we obtain
\begin{eqnarray*}
\int_{B_{1+\frac{\rho}{3}}(0)}|\nabla v_{n}|^{2}dx	 
&\leq& \displaystyle\int_{B_{2}(0)}|\nabla v_n|^{2}dx \\
&\leq& \displaystyle\int_{\bigcup_{k=1}^{m}B_{1}(\zeta_k)}|\nabla v_n|^{2}dx\\
&\leq& \displaystyle\sum_{k=1}^{m}\int_{B_{1}(\zeta_k)}|\nabla v_n|^{2}dx\\
&\leq& m\displaystyle\sup_{y\in\mathbb{R}^{N}}\int_{B_{1}(y)}|\nabla v_n|^{2}dx  
   \leq \frac{S_{H,L}^{(2N-\mu)/(N+2-\mu)}}{k}.
\end{eqnarray*}
Then,
\begin{eqnarray}\label{LP17}
\int_{B_{1+\frac{\rho}{3}}(0)}|\nabla v_{n}|^{2}dx \leq  \frac{S_{H,L}^{(2N-\mu)/(N+2-\mu)}}{k}, 
\end{eqnarray}
implying that 
$$
A \leq \frac{S_{H,L}^{(2N-\mu)/(N+2-\mu)}}{k}.
$$
Hence, by (\ref{m}),
$$
A < \left(\frac{S}{2^{\frac{2}{N}}}\right)^{\frac{2N-\mu}{N+2-\mu}}\left(\frac{1}{C(N,\mu)}\right)^{\frac{N-2}{N+2-\mu}},
$$
which contradicts (\ref{A}), and so, $v_0\neq 0$. 

Finally, using the equalities below, 
$$
\ds\int_{\mathbb{R}^{N}}|\nabla v_n|^{2}dx=\ds\int_{\mathbb{R}^{N}}|\nabla v_0|^{2}dx+\ds\int_{\mathbb{R}^{N}}|\nabla (v_n-v_0)|^{2}dx+o_{n}(1)
$$
and
$$
\ds\int_{\mathbb{R}^{N}}(I_{\mu}*|v_n|^{2^{*}_{\mu}})|v_n|^{2^{*}_{\mu}}dx=\ds\int_{\mathbb{R}^{N}}(I_{\mu}*|v_0|^{2^{*}_{\mu}})|v_0|^{2^{*}_{\mu}}dx+\ds\int_{\mathbb{R}^{N}}(I_{\mu}*|v_n-v_0|^{2{*}})|v_n-v_0|^{2^{*}}dx+o_{n}(1),
$$
(see \cite[Lemma 2.2]{Gao}) it follows that 
$$
E_{\infty}(v_n) =E_{\infty}(v_0)+E_{\infty}(v_n-v_0)+o_{n}(1).
$$
Since 
$$
E_\infty(u_n)=E_\infty(v_n),
$$
we have 
$$
E_{\infty}(u_n) =E_{\infty}(v_0)+E_{\infty}(v_n-v_0)+o_{n}(1),
$$
finishing the proof. 
{\fim}

The next result is crucial to study the compactness properties
involving the energy functional $E_\lambda$. We would like point out
that a version of that result for $\lambda=0$ can be found in
\cite[Lemma 3.1]{GSYZ} by using a different approach.

{\begin{teorema}\label{Teorema de Compacidade Global}(A global
    compactness result) Let  $(u_n)$ be a $(PS)_{c}$ sequence for
    $E_{\lambda}$ with $u_n \rightharpoonup u_0$ in $H$. Then, the
    sequence $(u_n)$ verifies either: \\ 
\noindent $(a)$ \,\, $u_n \to u_0$	\quad \mbox{or} \\	
\noindent $(b)$ there exist $k \in \mathbb{N}$ and $u^1,u^2,...,u^k$
nontrivial solution of $(P_\infty)$, such that 
\begin{eqnarray}\label{Conv12}
\|u_n\|^{2}=\|u_0\|^2+\displaystyle\sum^{k}_{j=1}\|u^{j}\|_{D^{1,2}}^{2}+o_{n}(1)
\end{eqnarray}
and
\begin{eqnarray}\label{Conv13}
E_{\lambda}(u_n)= E_{\lambda}(u_0) + \sum^{k}_{j=1}E_{\infty}(u^{j})+o_{n}(1).
\end{eqnarray}
\end{teorema}}

\noindent {\bf Proof.}  
Let us first consider the case $\lambda=0$. By Lemma \ref{PSbdd}, there exists $u_0 \in D^{1,2}(\R^N)$ such that
$u_n\rightharpoonup  u_0$ in $D^{1,2}$ and $u_{0}$  is a critical
point of  $E_{\lambda}$.  
Suppose that $u_n\nrightarrow u_{0}$ in $D^{1,2}(\mathbb{R}^{N})$ and
let  $(z^{1}_{n})\subset D^{1,2}(\mathbb{R}^{N})$ be the sequence
given by  $z^{1}_{n}=u_{n}-u_{0}$. 
Then,  $z^{1}_{n}\rightharpoonup 0$ in $D^{1,2}(\mathbb{R}^{N})$ and
$z^{1}_{n}\nrightarrow 0$ in  $D^{1,2}(\mathbb{R}^{N})$. 
Arguing as in the proof of Lemma \ref{Lematecnico}, we obtain  
\begin{eqnarray}\label{TCG7}
E_{\infty}(z^{1}_{n})=E_{0}(u_n)-E_{0}(u_0)+o_{n}(1)
\end{eqnarray}
and
\begin{eqnarray}\label{TCG12}
E'_{\infty}(z^{1}_{n})=E_{0}'(u_n)-E_{0}'(u_0)+o_{n}(1)\,\,\ {\rm in}\,\,\ (D^{1,2}(\mathbb{R}^{N}))'.
\end{eqnarray}
Then, from  $(\ref{TCG7})$ and  $(\ref{TCG12})$, we see that $(z^{1}_{n})$ is
a   $(PS)_{c_{1}}$ sequence for  $E_{\infty}$. Hence, by Lemma
$\ref{Lematecnico}$, there are sequences
$(\sigma_{n,1})\subset\mathbb{R}$, $(y_{n,1})\subset\mathbb{R}^{N}$ and  
a nontrivial solution $u^{1}\in D^{1,2}(\mathbb{R}^{N})$  for problem
$(P_{\infty})$ such that  
\begin{eqnarray}\label{TCG13}
v^{1}_{n}(x):=\sigma^{(N-2)/2}_{n,1}z^{1}_{n}\ds\bigg(\sigma_{n,1}x+y_{n,1}\ds\bigg) \rightharpoonup u^{1}_{0} \quad \mbox{in} \quad D^{1,2}(\mathbb{R}^N).
\end{eqnarray}
Since 
\begin{eqnarray}\label{TCG15}
\|v^{1}_{n}\|=\|u^{1}_{n}\|\,\,\ {\rm and}\,\,\ \displaystyle\int_{\mathbb{R}^{N}}(I_{\mu}*|v^{1}_{n}|^{2^{*}_{\mu}})|v^{1}_{n}|^{2^{*}_{\mu}} dx=\displaystyle\int_{\mathbb{R}^{N}}(I_{\mu}*|z^{1}_{n}|^{2^{*}_{\mu}})|z^{1}_{n}|^{2^{*}_{\mu}}dx,
\end{eqnarray}
it is easy to show that
\begin{eqnarray}\label{TCG16}
E_{\infty}(v^{1}_{n})=E_{\infty}(z^{1}_{n}) \quad \mbox{and} \quad E'_{\infty}(v^{1}_{n})\rightarrow 0\,\,\ {\rm in}\,\,\ (D^{1,2}(\mathbb{R}^{N}))',
\end{eqnarray}
hence $(v_n^1)$ is also a $(PS)_{c_1}$ for $E_\infty$.

Setting $z_n^{2}:=v_n^{1}-u^1$, we derive that
$$
E_\infty(v_n^{1})=E_\infty(z_n^{2})+E_\infty(u^{1})+o_n(1)
$$
$$
E'_\infty(v_n^{1})=E'_\infty(z_n^{2})+E'_\infty(u^{1})+o_n(1)\,\,\ {\rm in}\,\,\ (D^{1,2}(\mathbb{R}^{N}))',
$$
and
$$
\|v_n^1\|^{2}=\|z_n^2\|^{2}+\|u^1\|^{2}+o_n(1).
$$
Thus, 
$$
E_0(u_n)=E_0(u_0)+E_\infty(z_n^1)+o_n(1)=E_0(u_0)+E_\infty(v_n^1)+o_n(1)=E_0(u_0)+E_{\infty}(u^1)+E_\infty(z_n^2)+o_n(1)
$$
and
$$
\|u_n\|^{2}=\|u_0\|^{2}+\|z_n^1\|^{2}+o_n(1)=\|u_0\|^{2}+\|v_n^1\|^{2}+o_n(1)=\|u_0\|^{2}+\|u^1\|^{2}+\|z_n^{2}\|^{2}+o_n(1).
$$
Arguing as above, we have that $(z_n^{2})$ is a $(PS)_{c_2}$ sequence for $E_\infty$, that is, 
$$
E_\infty(z_n^{2}) \to c_2 \quad \mbox{and} \quad E'_\infty(z_n^{2}) \to 0. 
$$
Therefore, there are sequences
$(\sigma_{n,2})\subset\mathbb{R}$, $(y_{n,2})\subset\mathbb{R}^{N}$ and
 a nontrivial solution $u^{2}\in D^{1,2}(\mathbb{R}^{N})$  for problem
$(P_{\infty})$ such that  
\begin{eqnarray}\label{TCG13*}
v^{2}_{n}(x):=\sigma^{(N-2)/2}_{n,1}z^{1}_{n}\ds\bigg(\sigma_{n,1}x+y_{n,1}\ds\bigg) \rightharpoonup u^{2} \quad \mbox{in} \quad D^{1,2}(\mathbb{R}^N).
\end{eqnarray}

It is possible to prove that $(v_n^{2})$ is a $(PS)_{c_2}$ sequence
for $E_\infty$,  and fixing the sequence 
$$
z^{3}_{n}(x):=v^{2}_{n}(x)-u^{2}_{0}(x), 
$$
we will obtain  
$$
E_0(u_n)=E_0(u_0)+E_{\infty}(u^1)+E_{\infty}(u^2)+E_\infty(z_n^3)+o_n(1)
$$
and 
$$
\|u_n\|^{2}=\|u_0\|^{2}+\|u^1\|^{2}+\|u^2\|^{2}+\|z_n^{3}\|^{2}+o_n(1).
$$

Arguing as above, we will find  $u^{1}, u^{2},..., u^{k}$ nontrivial solutions for problem  $(P_{\infty})$ satisfying
\begin{eqnarray}\label{TCG34}
\|u_n\|^{2}=\|u_0\|^{2}+\ds\sum_{j=1}^{k}\|u^{j}\|^{2}+\|{z}^{k+1}_{n}\|^{2}+o_{n}(1)
\end{eqnarray}
and
\begin{eqnarray}\label{TCG35}
E_{0}(u_n)=E_{0}(u_0)+\ds\sum_{j=1}^{k}E_{\infty}(u^{j})+E_{\infty}({z}^{k+1}_{n})+o_{n}(1)
\end{eqnarray}
for the corresponding sequence $(z_n^{k+1})$ in $D^{1,2}(\mathbb{R}^N)$.

From the definition of $S_{H,L}$, 
\begin{eqnarray}\label{TCG36}
\biggl(\displaystyle\int_{\mathbb{R}^{N}}(I_{\mu}*|u^{j}|^{2^{*}_{\mu}})|u^{j}|^{2^{*}_{\mu}} dx\biggl)^{1/2^{*}_{\mu}}S_{H,L}\leq \|u^{j}\|^{2},\,\,\ j=1,2,...,k.
\end{eqnarray}
Since  $u^{j}$ is nontrivial solution of $(P_{\infty})$, for all $j=1,2,...,k$, we have
$$
\|u^{j}\|^{2}=\displaystyle\int_{\mathbb{R}^{N}}(I_{\mu}*|u^{j}|^{2^{*}_{\mu}})|u^{j}|^{2^{*}_{\mu}} dx.
$$
Hence, 
\begin{eqnarray}\label{TCG38}
\|u^{j}\|^{2} \geq  S_{H,L}^{(2N-\mu)/(N+2-\mu)},\,\,\ j=1,2,...,k.
\end{eqnarray}
From $(\ref{TCG34})$ and $(\ref{TCG38})$, 
\begin{eqnarray}\label{TCG39}
\|{z}^{k+1}_{n}\|^{2} = \|u_n\|^{2}-\|u_0\|^{2}-\ds\sum_{j=1}^{k}\|u^{j}\|^{2}+o_{n}(1)
\leq \|u_n\|^{2}-\|u_0\|^{2}-k\,S_{H,L}^{(2N-\mu)/(N+2-\mu)}+o_{n}(1).
\end{eqnarray}
Since  $(u_n)$ is bounded in  $D^{1,2}(\mathbb{R}^{N})$, for $k$
sufficient large, we conclude that  ${z}^{k+1}_{n}\rightarrow 0$
in $D^{1,2}(\mathbb{R}^{N})$, this proves the case $\lambda=0$. 

Let us now consider the case $\lambda>0$. As in the previous case, up to a subsequence $u_n\rightharpoonup  u_0$
in $H^1(\R^N)$, where $u_{0}\in H^1(\R^N)$   is a solution of $(P_\lambda)$.
We can also assume that  $u_n(x) \to u_0(x)$ a.e. in $\R^N$, and $u_n \to u_0$ in $L^2_{\loc}(\R^N)$.

Let us define $z_n^1:=u_n-u_0$.
Arguing as in the case $\lambda=0$, it follows that    
\beq
\label{1.2}
E_{\lambda,\infty}(z^1_{n})=E_{\lambda}(u_n)-E_{\lambda}(u_0)+o_{n}(1)
\eeq
and
\beq
\label{1.3}
E'_{\lambda,\infty}(z_{n}^1)=E_{\lambda}'(u_n)-E_{\lambda}'(u_0)+o_{n}(1)
\eeq
where
$$
E_{\lambda,\infty}(u):=E_\infty(u)+\frac{\lambda}{2}\int_{\R^N}u^2dx,\qquad u\in H^1(\R^N).
$$
Thus, it follows that $(z^{1}_{n})$ is a (PS) sequence for
$E_{\lambda,\infty}$ and
\beq
\label{1026}
z^{1}_{n} \to\quad{\rm weakly\  in }\,  H^{1}(\mathbb{R}^N),\ {\rm a.e.\ in }\ \R^N\ {\rm and\  in }\ L^2_{\loc}(\R^N). 
\eeq
If $z^1_n\to 0$ strongly in $H^1(\R^N)$, we are done, so let us assume
the existence of $c>0$ such that
\beq
\label{1.1}
\|z^1_n\|\ge c>0.
\eeq
Observe that (\ref{1.1}) implies the existence of a positive constant $\tilde c>0$ such that 
$$
|z^1_n|_{2^*}^{2^*}\ge \tilde c>0.
$$
Indeed, if this would be not the case, from  
\beq
\label{sbomb}
E'_{\lambda,\infty}(z^1_n)z^1_n=\int_{\R^N}(|\D z^1_n|^2+\lambda |z^1_n|^2)dx-\int_{\R^N} (I_{\mu}*|z^1_n|^{2^{*}_{\mu}})|z^1_n|^{2^{*}_{\mu}}dx=o_n(1),
\eeq
and (\ref{special}) we get at once $\|z^1_n\|\to 0$, as $n\to\infty$, contrary to (\ref{1.1}). 

Observe that
\beq
\label{nonvan}
 d^{1}_n:=\max_{i\in\NN}\|z^1_n\|_{L^{2^*}(Q_i)}>C>0,
 \eeq
where  $Q_i$, $i\in\NN$,  are hypercubes with disjoint interior and unitary sides such that $\R^N=\sum_{i\in \NN}Q_i$,
because 
 \begin{eqnarray*}
 0<\tilde c &\le& |z^1_n|_{2^*}^{2^*}=\sum\limits_{i=1}^\infty
 |z^1_n|_{L^{2^*}(Q_i)}^{2^*}\le (d_n^{1})^{2^*-2}\sum\limits_{i=1}^\infty
 |z^1_n|_{L^{2^*}(Q_i)}^{2}\\
 & \le & (d_n^{1})^{2^*-2} c_1 \sum\limits_{i=1}^\infty
 \|z^1_n\|_{H_1(Q_i)}^{2}\le c_2 (d_n^{1})^{2^*-2}.
 \end{eqnarray*}
 For every $n\in\NN$, let $y_n^{1}$ be the center of an hypercube where
 $d_n$ is attained and define
 \beq
 \label{305}
 \tilde{z}^1_n:=z^1_n(\cdot+y_n^{1}).
 \eeq 
It turns out that $(\tilde{z}^1_n)$ is a $(PS)$ sequence for $E_{\lambda,\infty}$, bounded in $H^1(\R^N)$.
Let us call $z$ the weak limit in $H^1(\R^N)$ of $(\tilde{z}^1_n)$, up to a subsequence.
Since $|y_n|\to\infty$ by (\ref{1026}), $z$  solves $(P_{\lambda,\infty})$, so that by Theorem \ref{nonET}
we get $z=0$.

Now we can argue exactly as in \cite[Lemma 3.3]{StruweBook}, and, taking into account of 
$|\tilde z^1_n|_{L^{2^*}(Q_0)}>C>0$ and of the other information we get,  we find a bounded sequence
of points $(x_n)$ in $\R^N$, an infinitesimal sequence $(\sigma_n)$ in $(0,+\infty)$ and a nontrivial solution 
$u^1$ of $(P_\infty)$ such that, if we define the sequence $(z^2_n)$ in $H^1(\R^N)$  by
\beq
\label{1059}
z^2_n(x):= \tilde z^1_n(x)-\varphi
\left(\frac{x-x_n}{\sigma_n^{1/2}}\right)\,
\frac{1}{\sigma^{\frac{N-2}{2}}}\, u^1\!\left(\frac{x-x_n}{\sigma_n}\right)\rightharpoonup 0\qquad{\rm in }\ H^1(\R^N),
\eeq
where $\varphi\in C^\infty_0(\R^N,[0,1])$ is a cut-off function such that $\varphi\equiv 1$ on $B_1(0)$, 
then $(z^2_n)$  is a (PS) sequence for $E_{\lambda,\infty}$, in $H^1(\R^N)$, that verifies
$$
z^2_n \rightharpoonup 0\qquad{\rm in }\ H^1(\R^N),
$$
and
\beq
\label{1100}
E_{\lambda,\infty}(z^2_n)=E_{\lambda,\infty}(\tilde z^1_n)-E_\infty(u^1)+o_n(1).
\eeq
We observe that in order to get (\ref{1100}) it is crucial that $\sigma_n\to 0$, which implies
\beq
\label{1116}
\int_{\R^N}(z^2_n)^2dx=\int_{\R^N}(\tilde z^1_n)^2dx+o(1)
\eeq
by (\ref{1059}). 
By (\ref{1116}) we can also write
\beq
\label{1216}
\|u_n\|^2_{H^1}=\|u_0\|^2_{H^1}+\|z^1_n\|_{H^1}^2+o_n(1)=
\|u_0\|^2_{H^1}+\|u^1\|_{D^{1,2}}^2+\|z^2_n\|_{H^1}^2+o_n(1).
\eeq
By (\ref{1.2}) and (\ref{1100}) we have
\beq
\label{1215}
E_\lambda(u_n)=E_\lambda(u_0)+E_\infty(u^1)+E_{\lambda,\infty}(z^2_n)+o_n(1).
\eeq
Moreover, since $(z^2_n)$ is a (PS) sequence for $E_{\lambda,\infty}$, it follows that 
$E_{\lambda,\infty}(z^2_n)\ge o_n(1)$, so that
$$
E_\lambda(u_n)\ge E_\lambda(u_0)+E_\infty (u^1)+o_n(1).
$$

If $z_n^2\to 0$ in $H^1(\R^N)$ we are done, otherwise we can iterate the procedure.
Taking into account that at every step $k$ we get
$$
E_\lambda(u_n)\ge E_{\lambda,\infty}(u_0)+\sum_{j=1}^k E_\infty(u^j)+o_n(1)
\ge E_{\lambda,\infty}(u_0)+k\, \frac{(N+2-\mu)}{4N-2\mu}\, S_{H,L}^{(2N-\mu)/(N+2-\mu)}+o_n(1),
$$
after a finite number of steps we reach a sequence $(z^{k+1}_n)$ such that 
$z^{k+1}_n\to 0$ in $H^1(\R^N)$.
Hence, we obtain (\ref{Conv12}) and (\ref{Conv13}) by iterating (\ref{1216}) and (\ref{1215}),
that completes the proof.

{\fim}

An immediate consequence of the last theorem are the next two corollaries.

\begin{corolario}\label{COR-1}
Let  $(u_n)$ be a  $(PS)_{c}$ sequence for  $E_{\lambda}$ with $c\in
\left(0,\displaystyle
\frac{(N+2-\mu)}{4N-2\mu}S_{H,L}^{(2N-\mu)/(N+2-\mu)}
\right)$. Then,
up to a subsequence,  $(u_n)$ strongly converges in  $H$. 
\end{corolario}

\noindent {\bf Proof.} Since $(u_n)$ is a (PS) sequence, $(u_n)$ is
bounded in  $H$, and so, for some subsequence, it follows that
$u_n\rightharpoonup u_0$ in $H$  
and $E'_{\lambda}(u_0)=0$ for some $u_0 \in H$. Suppose, by contradiction, that
$$
u_n\nrightarrow u_0\,\,\,\ {\rm in}\,\,\,\ H.
$$
From Theorem  $\ref{Teorema de Compacidade Global}$, there are
$k\in\mathbb{N}$ and nontrivial solutions  $z^{1}_{0}, z^{2}_{0},...,
z^{k}_{0}$ of problem $(P_{\infty})$ such that,  
$$
\|u_n\|^{2}\rightarrow \|u_0\|^{2}+\ds\sum_{j=1}^{k}\|z^{j}_{0}\|^{2}
$$
and
$$
E_{\lambda}(u_n)\rightarrow E_{\lambda}(u_0)+\ds\sum_{j=1}^{k}E_{\infty}(z^{j}_{0}).
$$
Note that
\begin{eqnarray*}
E_{\lambda}(u_0) &=&
\ds\frac{1}{2}\|u_0\|^{2}+\ds\frac{1}{2}\ds
\int_{\mathbb{R}^N}V_{\lambda}(x)u_0^{2}dx-\ds\frac{1}{2\cdot
  2^{*}_{\mu}}\displaystyle\int_{\mathbb{R}^{N}}
(I_{\mu}*|u_0|^{2^{*}_{\mu}})|u_0|^{2^{*}_{\mu}}dx
\\ &=& 
\ds\bigg(\ds\frac{1}{2}-\ds\frac{1}{2\cdot 2^{*}_{\mu}}\ds\bigg)
\displaystyle\int_{\mathbb{R}^{N}}
(I_{\mu}*|u_0|^{2^{*}_{\mu}})|u_0|^{2^{*}_{\mu}}dx\geq 0.
\end{eqnarray*}
Then, 
\begin{eqnarray*}
c & = & E_{\lambda}(u_0)+\ds\sum_{j=1}^{k}E_{\infty}(z^{j}_{0})\geq
\ds\sum_{j=1}^{k}E_{\infty}(z^{j}_{0})\\
& \geq &
k\ds\displaystyle\frac{(N+2-\mu)}{4N-2\mu}S_{H,L}^{(2N-\mu)/(N+2-\mu)}\geq
\displaystyle\frac{(N+2-\mu)}{4N-2\mu}S_{H,L}^{(2N-\mu)/(N+2-\mu)}, 
\end{eqnarray*}
which is a contradiction with  $c\in
\left(0,\displaystyle\frac{N+2-\mu}{4N-2\mu}S_{H,L}^{(2N-\mu)/(N+2-\mu)}\right)$. 
{\fim}


\begin{corolario}\label{COR-2}
The functional  $E_{\lambda}:H \rightarrow\mathbb{R}$ satisfies the
Palais-Smale condition in the range
$\left(\displaystyle\frac{(N+2-\mu)}{4N-2\mu}S_{H,L}^{(2N-\mu)/(N+2-\mu)},
2^{\frac{4-\mu}{N+2-\mu}}
\displaystyle\frac{(N+2-\mu)}{4N-2\mu}S_{H,L}^{(2N-\mu)/(N+2-\mu)}\right)$.  
\end{corolario}

\noindent {\bf Proof.} Let  $(u_n)$ be a  sequence in  $H$ that satisfies
$$
E_{\lambda}(u_n)\rightarrow c\,\,\,\ {\rm and}\,\,\,\
E_{\lambda}'(u_n)\rightarrow 0. 
$$
Since $(u_n)$ is bounded, up to a subsequence, we have
$u_n\rightharpoonup u_0$ in $H$, moreover $E_{\lambda}(u_0)\geq 0$. 
Suppose by contradiction that 
$$
u_n\nrightarrow u_0\,\,\,\ {\rm in}\,\,\,\ D^{1,2}(\mathbb{R}^N).
$$
From Theorem $\ref{Teorema de Compacidade Global}$, there are
$k\in\mathbb{N}$ and nontrivial solutions  $z^{1}_{0}, z^{2}_{0},...,
z^{k}_{0}$ of problem $(P_{\infty})$ such that  
$$
\|u_n\|^{2}\rightarrow \|u_0\|^{2}+\ds\sum_{j=1}^{k}\|z^{j}_{0}\|^{2}
$$
and
$$
E_{\lambda}(u_n)\rightarrow c=E_{\lambda}(u_0)+\ds\sum_{j=1}^{k}E_{\infty}(z^{j}_{0}).
$$
The above information ensures that $u_0 \not= 0$. Since
$E_{\lambda}(u_0)\geq 0$, then $k=1$ and  $z^{1}_{0}$ cannot change of
sign, because otherwise, by Lemma \ref{Nodal},   
$$
E_{\infty}(z^{j}_{0}) \geq 2^{\frac{4-\mu}{N+2-\mu}}
\displaystyle\frac{(N+2-\mu)}{4N-2\mu}S_{H,L}^{(2N-\mu)/(N+2-\mu)}, 
$$
which leads to a contradiction. Thereby, as $z^{1}_{0}$ has definite
sign, $z^{1}_{0}=U_{\sigma,z}$ for suitable $\sigma>0$ and $z\in\R^N$
and, by (\ref{minimodeI_0}),
$$
E_{\infty}(z^{j}_{0}) =
\displaystyle\frac{(N+2-\mu)}{4N-2\mu}S_{H,L}^{(2N-\mu)/(N+2-\mu)}. 
$$
On the other hand, by a direct computation, 
$$
E_{\lambda}(u_{0}) \geq
\displaystyle\frac{(N+2-\mu)}{4N-2\mu}S_{H,L}^{(2N-\mu)/(N+2-\mu)}. 
$$
Hence, 
$$
c=E_{\lambda}(u_0)+E_{\infty}(z^{1}_{0}) \geq
2\displaystyle\frac{(N+2-\mu)}{4N-2\mu}S_{H,L}^{(2N-\mu)/(N+2-\mu)}
>
2^{\frac{4-\mu}{N+2-\mu}}
\displaystyle\frac{(N+2-\mu)}{4N-2\mu}S_{H,L}^{(2N-\mu)/(N+2-\mu)}, 
$$
obtaining again a contradiction. This proves the result.
{\fim}

\vspace{2mm}


The next results provide us the $(PS)$ condition for the functional
$J_\lambda$.  
The first one is a direct computation and we omit its
proof, the second one is an immediate consequence of the study made above. 

\begin{lema}\label{LLL-1}
Let $(u_n)\subset\mathcal{M}$ be a sequence that satisfies 
$$
J_{\lambda}(u_n)\rightarrow c\,\,\,\ and \,\,\,\
J_{\lambda}'|_{\mathcal{M}}(u_n)\rightarrow 0. 
$$
Then, the sequence  $v_n=c^{(N-2)/(2N-2\mu +4)}u_n$ satisfies the following limits. 
$$
E_{\lambda}(v_n)\rightarrow
\ds\frac{(N+2-\mu)}{4N-2\mu}c^{(2N-\mu)/(N+2-\mu)}\,\,\,\ and\,\,\,\
E_{\lambda}'(v_n)\rightarrow 0. 
$$
\end{lema}


\begin{corolario}\label{LLL-2}
Suppose that there are a sequence  $(u_n)\subset\mathcal{M}$ and 
$$
c\in (S_{H,L},2^{\frac{4-\mu}{2N-\mu}}S_{H,L})
$$
such that
$$
J_{\lambda}(u_n)\rightarrow c\,\,\,\ and\,\,\,\
J_{\lambda}'|_{\mathcal{M}}(u_n)\rightarrow 0. 
$$
Then
\begin{itemize}
\item[{\em a)}] there exists $u_0\in\cM$ such that, up to a
  subsequence, $u_n\rightarrow u_0$ in $D^{1,2}(\mathbb{R}^N)$ and
  $u_0$ is a critical point for $J_\lambda$ constrained on $\cM$;
\item[{\em b)}]
$E_{\lambda}$ has a critical point  $v_0\in H$ with 
$E_{\lambda}(v_0)=\displaystyle\frac{(N+2-\mu)}{4N-2\mu}c^{(2N-\mu)/(N+2-\mu)}$.
\end{itemize}
\end{corolario}



\section{Main tools and basic estimates }

We are looking for solutions of problem $(P_\lambda)$ as critical
points of the functional $J_\lambda$ constrained on $\cM$, up to a multiplier.
Next proposition shows that the problem cannot be solved by
minimization, so no ground state solution exists.

\begin{proposicao}\label{Proposition3.1}
Set 
\begin{eqnarray}\label{minimization}
m:=\inf\{J_{\lambda}(u): u\in \mathcal{M}\}.
\end{eqnarray} 
Then
$$
m=S_{H,L}
$$
and the  minimization problem (\ref{minimization}) has no solution.
\end{proposicao}
\noindent {\bf Proof}.  
Let $u \in \mathcal{M}$ be  arbitrarily chosen. Then, by $(V_1)$ we get
$$
J_{\lambda}(u) \geq S_{H,L},
$$
which implies
$$
m \geq S_{H,L}.
$$

In order to show the opposite inequality, let us consider the sequence
$$
\tilde{\Psi}_{n}(x) =\xi(|x|)U_{\frac{1}{n},0}(x),
$$
where $\xi \in C^{\infty}_{0}((0,+\infty),[0,1])$ is such that $\xi(s)=1$ is $s\in [0,1/2]$ and $\xi(s)=0$ is $s\geq 1$. Using (\ref{vouusar}) together with the definition of $\tilde{\Psi}_{n}$, we have 
\begin{eqnarray}\label{grad}
\displaystyle\int_{\mathbb{R}^{N}}|\nabla \tilde{\Psi}_{n}|^{2} dx=\displaystyle\int_{\mathbb{R}^{N}}|\nabla U_{\frac{1}{n},0}|^{2} dx + o_{n}(1),
\end{eqnarray}
\begin{eqnarray}\label{2*}
\int_{\mathbb{R}^N}(I_{\mu}*|\tilde{\Psi}_{n}|^{2^{*}_{\mu}})|\tilde{\Psi}_{n}|^{2^{*}_{\mu}}dx=\int_{\mathbb{R}^N}(I_{\mu}*|U_{\frac{1}{n},0}|^{2^{*}_{\mu}})U_{\frac{1}{n},0}|^{2^{*}_{\mu}}dx + o_{n}(1)
\end{eqnarray}
and 
\begin{eqnarray}\label{L2}
\lambda\displaystyle\int_{\mathbb{R}^{N}}| \tilde{\Psi}_{n}|^{2} dx= o_{n}(1). 
\end{eqnarray}
On the other hand, for all $\rho>0$, we have
\begin{eqnarray*}
	\displaystyle\int_{\mathbb{R}^{N}}V_{0}(x)|\tilde{\Psi}_{n}|^{2} dx &=& \displaystyle\int_{B_{\rho}(0)}V_{0}(x)|\tilde{\Psi}_{n}|^{2} dx+\displaystyle\int_{\mathbb{R}^{N}\setminus B_{\rho}(0)}V_{0}(x)|\tilde{\Psi}_{n}|^{2} dx\\
	&\leq& | \tilde{\Psi}_{n}|^{2}_{L^{2^{*}}(\mathbb{R}^{N})}\displaystyle\left(\displaystyle\int_{B_{\rho}(0)}
	|V_{0}(x)|^{N/2}dx
	\displaystyle\right)^{2/N}\\
	&+& |V_{0}|_{L^{N/2}(\mathbb{R}^{N})}\displaystyle\left(\displaystyle\int_{\mathbb{R}^{N}\setminus B_{\rho}(0)}| \tilde{\Psi}_{n}(x)|^{2^{*}}dx\ds\right)^{2/2^{*}}.
\end{eqnarray*}
Now, recalling that 
$$
\ds\lim_{n\rightarrow \infty}\int_{\mathbb{R}^{N}\setminus B_{\rho}(0)}| \tilde{\Psi}_{n}(x)|^{2^{*}}dx=0,
$$
$$
\sup_{n \in \mathbb{N}}| \tilde{\Psi}_{n}|_{2^{*}}<+\infty
$$
and
$$
\ds\lim_{\rho \to 0}\displaystyle\int_{B_{\rho}(0)}
|V_{0}(x)|^{N/2}dx
=0,
$$
we get
\begin{eqnarray}\label{PROP1-P2}
\ds\lim_{n\rightarrow \infty}\ds\int_{\mathbb{R}^{N}}V_{0}(x) |\tilde{\Psi}_{n}(x)|^{2}dx=0.
\end{eqnarray}
Now, if we define
$$
\hat\Psi_n(x)=\frac{1}{\int_{\mathbb{R}^N}
(I_{\mu}*|\tilde{\Psi}_{n}|^{2^{*}_{\mu}})|\tilde{\Psi}_{n}|^{2^{*}_{\mu}}dx}\, 
\tilde\Psi_n(x) 
$$
then $\hat\Psi_n\in\mathcal{M}$, $\forall n\in\N$, and from
$(\ref{grad})-(\ref{PROP1-P2})$,  
$$
\ds\lim_{n\rightarrow \infty}J_{\lambda}(\hat{\Psi}_{n}(x))=S_{H,L},
$$
which concludes the first part of the  proof. Now suppose that the
minimization problem $(\ref{minimization})$ has a solution $u^{*}$. 
Then
$$
S_{H,L}\leq \ds\frac{\ds\int_{\mathbb{R}^{N}}|\nabla
  u^{*}|^{2}dx}{\biggl(\displaystyle
\int_{\mathbb{R}^{N}}(I_{\mu}*|u|^{2^{*}_{\mu}})|u|^{2^{*}_{\mu}}
  dx\biggl)^{1/2^{*}_{\mu}}}\leq
\ds\frac{\ds\int_{\mathbb{R}^{N}}[|\nabla
  u^{*}|^{2}+V_{\lambda}(x)|u^{*}|^{2}]dx}{\biggl(
\displaystyle\int_{\mathbb{R}^{N}}(I_{\mu}*|u^{*}|^{2^{*}_{\mu}})|u^{*}|^{2^{*}_{\mu}}
  dx\biggl)^{1/2^{*}_{\mu}}}=S_{H,L}. 
$$
The above relation implies that
$\ds\int_{\mathbb{R}^{N}}V_{\lambda}(x)|u^{*}|^{2}dx=0$ and
$u^{*}=U_{\sigma,z}$ 
for some $\sigma>0$ and $z \in \mathbb{R}^N$. Thus, using the
assumptions on $V_\lambda$ and the fact that $U_{\sigma,z}>0$ for all
$x\in\mathbb{R}^{N}$, we deduce 
$$
0=\ds\int_{\mathbb{R}^{N}}V_{\lambda}(x)|u^{*}|^{2}dx=
\ds\int_{\mathbb{R}^{N}}V_{\lambda}(x)|U_{\sigma,z}|^{2}dx>0, 
$$
which is impossible.

{\fim}

In view of the previous proposition, the main goal of this section
will be to introduce some tools and to establish some basic estimates
in oder to find bound state solutions in the next section. 
To begin with, let us introduce a barycenter type map $\beta:
H\setminus\{0\}\to \mathbb{R}^{ N}$ given by 
$$\beta(u)=\displaystyle\frac{1}{|u|_{2^*}^{2^{*}}}\ds
\displaystyle\int_{\mathbb{R}^{N}}\frac{x}{1+|x|}|u|^{2^{*}}dx 
$$
and a kind of inertial momentum $\gamma: H\setminus\{0\}\to \mathbb{R}$ given by
$$\gamma(u)=\displaystyle \frac{1}{|u|_{2^*}^{2^{*}}}\ds\displaystyle\int_{\mathbb{R}^{N}}\biggl|\frac{x}{1+|x|}-\beta(u)\biggl||u|^{2^{*}}dx.
$$
It is readily seen that the maps $\beta$ and $\gamma$ are continuous and, moreover, 
$\beta(tu)=\beta(u)$ and $\gamma(tu)=\gamma(u)$, for all $t\in \mathbb{R}$ and for all $u\in H\setminus\{0\}$.

\begin{lema}\label{Proposition4.1}
Let $\lambda\ge 0$ and define
$$
\mathcal{B}_{V_{\lambda}}=\inf\left\{J_{\lambda}(u): u \in \mathcal{M}, \beta(u)=0, \gamma(u)=\frac 12\right\}.
$$
Then the following inequalities hold
$$
\mathcal{B}_{V_{\lambda}} > S_{H,L}, \ \ \mbox{for all} \ \ \lambda \geq 0.
$$
\end{lema}
\noindent {\bf Proof.} By Proposition \ref{Proposition3.1}, 
$$
\inf\left\{J_{\lambda}(u): u \in \mathcal{M}, \beta(u)=0, \gamma(u)=\frac 12\right\}\geq S_{H,L}.
$$
Now suppose, by contradiction, that the equality is true. Then, there exists a sequence $(u_n)\subset H$ such that 
\begin{equation}\label{PROP3-P1}
\begin{cases}
(a) \ \ u_n \in \mathcal{M},\ \ \beta(u_n)=0,\ \ \gamma(u_n)=\ds\frac{1}{2};\\
(b) \ \ \ds\lim_{n\rightarrow+\infty}J_{\lambda}(u_n)=S_{H,L}.
\end{cases}
\end{equation}
Note that
\begin{eqnarray*}
 S_{H,L}= J_{\lambda}(u_n)+o_{n}(1) \geq
 \displaystyle\int_{\mathbb{R}^{N}}|\nabla u_n|^{2} dx +o_n(1)\geq
 S_{H,L}+o_n(1), 
 \end{eqnarray*}
then
\begin{eqnarray} \label{NEWE1} 
\displaystyle\int_{\mathbb{R}^{N}}|\nabla u_n|^{2} dx +o_n(1)=  S_{H,L}.
\end{eqnarray}

By Lemma \ref{SHLS} and \cite[Theorem 2.5]{BC}, we get
\begin{eqnarray}\label{fromunicity}
u_n(x)=\Theta \, {U_{\sigma_n,y_n}(x)}+ \epsilon_{n}(x)
\end{eqnarray}
for some positive constant $\Theta$ with $U_{\sigma_n,y_n}$ as in (\ref{U}), $\sigma_n>0$, $y_n \in
\mathbb{R}^{N}$ and  $\epsilon_{n} \to 0$ in
$D^{1,2}(\mathbb{R}^{N})$. 
In order to get the constant $\Theta>0$, we recall that $|U_{\sigma_n,y_n}|_{2^*}=|U_{1,0}|_{2^*}$ and that by (\ref{NEWE1})  the sequence $(u_n)$ is bounded
in $D^{1,2}(\mathbb{R}^N)$, and so, we can assume that for some subsequence  $\displaystyle \lim_{n \to +\infty}|u_n|_{2^*}=L_1$. Moreover, since $u_n \in \mathcal{M}$, we must have $L_1>0$. 

We claim that $\displaystyle\lim_{n\to \infty}\sigma_n
=\overline{\sigma}>0$ and $\displaystyle\lim_{n\to \infty}y_n
=\overline{y}$ in $\mathbb{R}^{N}$. Let us first show that
$(\sigma_n)$ is bounded. In fact, if for some subsequence, still
denoted by $(\sigma_n)$,
$\ds\lim_{n\rightarrow+\infty}\sigma_n=+\infty$ occurs, then  for all
$\rho>0$, we have 
\begin{eqnarray*}
\ds\lim_{n\rightarrow+\infty}\int_{B_{\rho}(0)}|u_n|^{2^{*}}dx 
&=& 
\Theta^{2^*}\, \lim_{n\rightarrow+\infty}\int_{B_{\rho}(0)}|U_{\sigma_n,y_n}|^{2^{*}}dx=0.
\end{eqnarray*}
Since $\beta(u_n)=0$, for all $\rho>0$, 
\begin{eqnarray*}
\gamma(u_n) 
&=& 
\frac{1}{|u_n|_{2^*}^{2^*}}
\displaystyle\int_{\mathbb{R}^{N}}\frac{|x|}{1+|x|}|u_{n}|^{2^{*}}dx
\\
&= &  
\frac{1}{|u_n|_{2^*}^{2^*}}
\displaystyle\int_{\mathbb{R}^{N}\setminus
  B_{\rho}(0)}\frac{|x|}{1+|x|}|u_{n}|^{2^{*}}dx
\\ 
& &  +
\frac{1}{|u_n|_{2^*}^{2^*}}
\displaystyle\int_{B_{\rho}(0)}\frac{|x|}{1+|x|}|u_{n}|^{2^{*}}dx
\\
&\geq& 
\ds\frac{\rho}{1+\rho}+o_{n}(1),
\end{eqnarray*}
so
$$
\ds\liminf_{n\rightarrow+\infty}\gamma(u_n)\geq \ds\frac{\rho}{1+\rho}, \ \ \forall \rho>0,
$$
and then
\begin{eqnarray}\label{LEMA2-P3}
\ds\liminf_{n\rightarrow+\infty}\gamma(u_n)\geq 1,
\end{eqnarray}
obtaining, therefore, a contradiction. Thus, $(\sigma_n)$ is bounded
and we can assume that 
$$
\ds\lim_{n\rightarrow+\infty}\sigma_n=\overline{\sigma} \ \ {\rm with}
\ \ \overline{\sigma}\geq 0. 
$$
We claim that $\overline{\sigma}$ is positive. In fact, if
$\overline{\sigma}=0$,  for all $\rho>0$ we have  
\begin{eqnarray*}
\ds\lim_{n\rightarrow+\infty}\int_{\mathbb{R}^{N}\setminus
  B_{\rho}(y_n)}|u_n|^{2^{*}}dx 
&=&
\Theta^{2^*}\, \lim_{n\rightarrow+\infty}\int_{\mathbb{R}^{N}\setminus
  B_{\rho}(y_n)}|U_{\sigma_n,y_n}|^{2^{*}}dx=0. 
\end{eqnarray*}
As $\beta(u_n)=0$, we get
\begin{eqnarray}
\nonumber \ds\frac{|y_n|}{1+|y_n|} 
&=& 
\biggl|\ds\frac{y_n}{1+|y_n|}
-\beta(u_n)\biggr|=\frac{1}{|u_n|_{2^*}^{2^*}}
\ds\left|\ds\int_{\mathbb{R}^{N}}\ds\left(
\frac{y_n}{1+|y_n|}-\frac{x}{1+|x|}\ds\right)|u_n|^{2^{*}}dx\ds\right|
\nonumber\\ 
\label{1132} &\leq& 
\frac{1}{|u_n|_{2^*}^{2^*}}
\ds\left|\ds\int_{B_{\rho}(y_n)}
\ds\left(\frac{y_n}{1+|y_n|}-\frac{x}{1+|x|}\ds\right)|u_n|^{2^{*}}dx\ds\right|
\\ 
\nonumber & & +
\frac{1}{|u_n|_{2^*}^{2^*}}
\ds\left|\ds\int_{\mathbb{R}^{N}\setminus
    B_{\rho}(y_n)}\ds
\left(\frac{y_n}{1+|y_n|}-\frac{x}{1+|x|}\ds\right)|u_n|^{2^{*}}dx\ds\right|
\leq 
\rho + o_n(1).
\end{eqnarray}
Hence,
$$
\ds\limsup_{n\rightarrow+\infty}\ds\frac{|y_n|}{1+|y_n|}\leq 
\rho, \ \ \forall \rho>0,
$$
from which it follows
$$
\ds\lim_{n\rightarrow+\infty}|y_n|=0.
$$
On the other hand, by the same calculus performed in (\ref{1132}),
\begin{eqnarray*}
\ds\lim_{n\rightarrow+\infty}\gamma(u_n) &=&
\ds\lim_{n\rightarrow+\infty}\frac{1}{|u_n|_{2^*}^{2^*}}
\int_{\mathbb{R}^{N}}\ds\left|
\frac{x}{1+|x|}-\beta(u_n)\ds\right||u_n|^{2^{*}}dx\nonumber\\ 
&=& \ds\lim_{n\rightarrow+\infty}\frac{1}{|u_n|_{2^*}^{2^*}}
\int_{\mathbb{R}^{N}}\ds\left|\frac{x}{1+|x|}-\frac{y_n}{1+|y_n|}\ds
\right||u_n|^{2^{*}}dx=0, 
\end{eqnarray*}
which is a contradiction.

Now, we are able to prove that $(y_n)$ is bounded. For this, suppose
by contradiction, that there is a subsequence, still denoted by
$(y_n)$, such that 
$$
\ds\lim_{n\rightarrow+\infty}|y_n|=+\infty.
$$
Then, for all $\epsilon>0$, there is $R>0$ and $n_0\in\mathbb{N}$ such that
\begin{eqnarray}\label{LEMA2-P8}
|x-y_n|<R \Rightarrow
\ds\left|\frac{x}{1+|x|}-\frac{y_n}{1+|y_n|}\ds\right|<\epsilon, \ \
\forall n\geq n_0 
\end{eqnarray}
and
\begin{eqnarray}\label{LEMA2-P9}
\ds\int_{\mathbb{R}^{N}\setminus
  B_{R}(y_n)}|u_{n}|^{2^{*}}dx=\Theta^{2^*}\, \int_{\mathbb{R}^{N}\setminus
  B_{R}(y_n)}|U_{\sigma_n,y_n}|^{2^{*}}dx+o_n(1)<\epsilon. 
\end{eqnarray}
From (\ref{LEMA2-P8}) and (\ref{LEMA2-P9}),
\begin{eqnarray*}
\ds\left|\beta(u_n)-\frac{y_n}{1+|y_n|}\ds\right| 
&\leq& 
\ds\frac{1}{|u_n|_{2^*}^{2^*}}
\int_{\mathbb{R}^{N}}\ds\left|\frac{x}{1+|x|}-
\frac{y_n}{1+|y_n|}\ds\right||u_n|^{2^{*}}dx
\\
&=& 
\frac{1}{|u_n|_{2^*}^{2^*}}
\ds\int_{B_{R}(y_n)}\ds\left|\frac{x}{1+|x|}
-\frac{y_n}{1+|y_n|}\ds\right||u_n|^{2^{*}}dx 
\\
&& + 
\frac{1}{|u_n|_{2^*}^{2^*}} \ds\int_{\mathbb{R}^{N}\setminus
  B_{R}(y_n)}\ds\left|\frac{x}{1+|x|}-
\frac{y_n}{1+|y_n|}\ds\right||u_{n}|^{2^{*}}dx
\\
&\leq& 
\epsilon+2\, \frac{\epsilon}{|u_n|_{2^*}^{2^*}}+o_n(1)
\end{eqnarray*}
which implies
$$
\ds\lim_{n\rightarrow+\infty}|\beta(u_n)|=1,
$$
which again leads us to a contradiction. 
Therefore, $(y_n)$ is bounded and we can assume that
$$
\ds\lim_{n\rightarrow+\infty}y_n=\overline{y}.
$$
Then
\begin{eqnarray*}
S_{H,L}
&=&
\displaystyle\lim_{n\to\infty}\left[\int_{\mathbb{R}^{N}}(|\nabla
u_n|^{2}+ V_{\lambda}(x)u_{n}^{2})dx\right]
\\
&=&
\Theta^{2}\, \left[\int_{\mathbb{R}^{N}}(|\nabla
U_{\overline{\sigma},\overline{y}}|^{2}+
V_{\lambda}(x)|U_{\overline{\sigma},\overline{y}}|^{2})dx\right]
\\
&>& 
\Theta^{2}\, 
\int_{\mathbb{R}^{N}}|\nabla U_{\overline{\sigma},\overline{y}}|^{2}dx =S_{H,L},
\end{eqnarray*}
which is an absurd. 
{\fim}


\begin{lema}\label{Proposition3.2}
If $\lambda>0$,  then
\beq
\label{1145}
\Upsilon_\lambda:=\inf\{J_{\lambda}(u): u \in \mathcal{M}, \beta(u)=0, 
\gamma(u)\geq 1/2\}> S_{H,L}.
\eeq
\end{lema}
\noindent {\bf Proof.} We start observing that  
$$
\Upsilon_\lambda=\inf\{J_{\lambda}(u): u \in \mathcal{M}, \beta(u)=0,
\gamma(u)\geq 1/2\}\geq S_{H,L}. 
$$
Now suppose, by contradiction, that the equality is true. Then, there
exists a sequence $(u_n)$ such that  
\begin{equation}
\begin{cases}
(a) \ \ u_n \in \mathcal{M},\ \ \beta(u_n)=0,\ \ \gamma(u_n)\geq\ds {1}/{2};\\
(b) \ \ \ds\lim_{n\rightarrow+\infty}J_{\lambda}(u_n)= S_{H,L}.
\end{cases}
\end{equation}

Then, the same computations made in Lemma \ref{Proposition4.1} allow
to assert that  
\beq
u_n(x)=\Theta \, U_{\sigma_n,y_n}(x)+ \epsilon_{n}(x)
\eeq
with $\Theta>0$, $\sigma_n>0$, $y_n \in \mathbb{R}^{N}$,
$\epsilon_{n} \to 0$ in $D^{1,2}(\mathbb{R}^{N})$ verifying
$\displaystyle\lim_{n\to \infty}\sigma_n
=\overline{\sigma}\in(0,+\infty]$ and $\displaystyle\lim_{n\to
  \infty}y_n =\overline{y}$ in $\mathbb{R}^{N}$. 
Let us show that $\overline\sigma=+\infty$ cannot occur.
If this would be the case, then
\begin{eqnarray*}
S_{H,L} 
&\ge & 
\liminf_{n\to+\infty}\left[\int_{\R^N}|\nabla
  u_n|^2dx+\lambda\int_{B_{\sqrt\sigma_n}(y_n)}|u_n|^2dx\right]
\\
&\ge &
\left[S_{H,L}+\lambda \Theta^{2}  \displaystyle \liminf_{n\to+\infty}\sigma_n \, \int_{B_1(0)}|U_{1,0}|^2dx\right]
\\
&=&
+\infty
\end{eqnarray*}
that is a contradiction. 
So, we can assume that $\overline\sigma\in(0,\infty)$, and then
\begin{eqnarray*}
S_{H,L}
&\geq&
\displaystyle\lim_{n\to\infty}\biggl[\int_{\mathbb{R}^{N}}(|\nabla
u_n|^{2}+ \lambda |u_{n}|^{2})dx\biggl]
\\
&\ge &
\Theta^2\biggl[\int_{\mathbb{R}^{N}}|\nabla
U_{\overline{\sigma},\overline{y}}|^{2}\,dx+ {\lambda} \overline\sigma
\int_{B_{\overline{\sigma}}(\overline{y})}|U_{\overline{\sigma},\overline{y}}|^{2}dx
 \biggr]
\\
&>& 
\Theta^2 \int_{\mathbb{R}^{N}}|\nabla
U_{\overline{\sigma},\overline{y}}|^{2}dx 
=S_{H,L},
\end{eqnarray*}
which is again a contradiction.  {\fim}

\begin{obs}
Testing the functional $J_0$ by the functions $\frac{1}{N}\, U_{n,0}$,
$n\in\NN$, it 
is readily seen that 
$$
\inf\left\{J_{0}(u): u \in \mathcal{M}, \beta(u)=0, \gamma(u) \geq
  1/2\right\}= S_{H,L}. 
$$ 
\end{obs}


\vspace{3mm}


Let $a\in (0,1)$ be such that 
\beq
\label{1930}
|V_0|_{L^{N/2}(\mathbb{R}^{N})}=
S\, \left(2^{a\, \frac{4-\mu}{2N-\mu}}-1\right) 
\eeq
and let us fix a number $\overline{c}$ such that 

\begin{equation} \label{cbarra}
S_{H,L}<\overline{c}< \min \biggl(
\frac{\mathcal{B}_{V_{0}}+ S_{H,L}}{2}, 
2^{(1-a)\, \frac{4-\mu}{2N-\mu}}S_{HL}
\biggr). 
\end{equation}
Note that this interval is not empty by Lemma \ref{Proposition4.1}.

In the sequel, $\varphi$ is a function that belongs to
$C_{0}^{\infty}(B_{1}(0))$ and satisfyies the following properties: 
\begin{equation}\label{PROPRIEDADE DA VARPHI}
\begin{cases}
\varphi\in C^{\infty}_{0}(B_{1}(0)),\ \ \varphi(x)>0\ \ \forall x\in B_{1}(0),\\
\varphi\ \ {\rm is\ \ symmetric\ \ and}\ \ |x_1|<|x_2|\Rightarrow
\varphi(x_1)>\varphi(x_2),
\\ 
\varphi\in \mathcal{M} \ \  \mbox{and} \ \ 
\displaystyle\int_{\mathbb{R}^{N}}|\nabla \varphi|^{2} dx=\Sigma \in
(S_{H,L}, \overline{c}). 
\end{cases}
\end{equation}

For every $\sigma>0$ and $y\in\mathbb{R}^{N}$, we set 
\begin{equation}\label{TRANSLACAO}
\varphi_{\sigma,y}(x)=
    \begin{cases}
    \sigma^{-\frac{N-2}{2}}\varphi\ds\left(\ds\frac{x-y}{\sigma}\ds\right),\
    \ x\in B_{\sigma}(y),\\ 
    0,\ \ \ \ \ \ \ \ \ \ \ \ \ \ \ \ \ \ \ \ \ \ \ \ \ x\notin B_{\sigma}(y).
    \end{cases}
\end{equation}
We remark that by the definition of $\varphi_{\sigma,y}$ and by
variable change, it follows that for every $\sigma>0$ and $y\in\R^N$
\begin{eqnarray}\label{familiadecatoff1}
\displaystyle\int_{\mathbb{R}^{N}}|\nabla \varphi_{\sigma,y}|^{2}
dx=\displaystyle\int_{B_{\sigma}(y)}|\nabla \varphi_{\sigma,y}|^{2}
dx=\displaystyle\int_{B_{1}(0)}|\nabla \varphi|^{2} dx, 
\end{eqnarray}
\begin{eqnarray}\label{familiadecatoff}
\displaystyle\int_{\mathbb{R}^{N}}|\varphi_{\sigma,y}|^{2^{*}}
dx=\displaystyle\int_{B_{\sigma}(y)}|\varphi_{\sigma,y}|^{2^{*}}
dx=\displaystyle\int_{B_{1}(0)}|\varphi|^{2^{*}} dx 
\end{eqnarray}
and
$$
\int_{\mathbb{R}^N}(I_{\mu}*|\varphi_{\sigma,y}|^{2^{*}_{\mu}})
\, |\varphi_{\sigma,y}|^{2^{*}_{\mu}}dx
=
\int_{\mathbb{R}^N}(I_{\mu}*\varphi^{2^{*}_{\mu}})
\, \varphi^{2^{*}_{\mu}}dx
=1
$$
so that, in particular,
\beq
\label{1826}
\varphi_{\sigma,y}\in\mathcal{M}\quad \mbox{ and }\quad
\int_{\R^N}|\nabla \varphi_{\sigma,y}|^2dx=\Sigma\in (S_{H,L},
\overline{c}) 
\qquad\forall\sigma>0\ \mbox{ and }\ \forall y\in\R^N.
\eeq


\begin{lema}\label{LEMA3}
The following relations hold:
\begin{equation}
    \begin{cases}
    (a)\ \ \ds\lim_{\sigma\rightarrow
      0}\sup\ds\left\{\ds\int_{\mathbb{R}^{N}}V_{0}(x)|\varphi_{\sigma,y}|^2dx;\
      \ y\in\mathbb{R}^{N}\ds\right\}=0;
\\
    (b)\ \
    \ds\lim_{\sigma\rightarrow+\infty}\sup\ds\left\{
    \ds\int_{\mathbb{R}^{N}}V_{0}(x)|\varphi_{\sigma,y}|^2\,dx;\
      \ y\in \mathbb{R}^{N}\ds\right\}=0;
\\
    (c)\ \
    \ds\lim_{r\rightarrow+\infty}\sup\ds\left\{
    \ds\int_{\mathbb{R}^{N}}V_{0}(x)|\varphi_{\sigma,y}|^2\,dx;\
      \ |y|=r,\ \ \sigma>0,\ \  y\in\mathbb{R}^{N}\ds\right\}=0. 
    \end{cases}
\end{equation}
\end{lema}
\noindent {\bf Proof.}  \, 
Let $y\in\mathbb{R}^{N}$ be chosen arbitrarily. Then, by the H\"older
inequality, 
\begin{eqnarray*}
    \ds\int_{\mathbb{R}^{N}}V_{0}(x)|\varphi_{\sigma,y}|^2\,dx &=&
    \ds\int_{ B_{\sigma}(y)}V_{0}(x)|\varphi_{\sigma,y}|^2\,dx\leq
    |V_{0}|_{L^{N/2}(
      B_{\sigma}(y))}|\varphi_{\sigma,y}|^{2}_{L^{2^{*}}(B_{\sigma}(y))}
\\
    &=& 
|V_{0}|_{L^{N/2}( B_{\sigma}(y))}|\varphi|^{2}_{2^{*}},\qquad \forall\sigma>0,
\end{eqnarray*}
hence
\begin{equation}\label{L3-P1}
\ds\sup\left\{\ds\int_{\mathbb{R}^{N}}V_{0}(x)|\varphi_{\sigma,y}|^2\,dx;\
  \ y\in\mathbb{R}^{N}\ds\right\}\leq
|\varphi|^{2}_{L^{2^{*}}(B_{1}(0))} \ds\sup \left\{|V_{0}|_{L^{N/2}(
    B_{\sigma}(y))};\ \ y\in\mathbb{R}^{N}\ds\right\}. 
\end{equation}
Since
$$
\ds\lim_{\sigma\rightarrow 0}\sup_{y \in \mathbb{R}^N}|V_{0}|_{L^{N/2}( B_{\sigma}(y))}=0,
$$
so $(a)$ follows from $(\ref{L3-P1})$.

To prove $(b)$, we fix arbitrarily $y\in\mathbb{R}^{N}$ and note that
by the H\"older inequality,  
\begin{eqnarray*}
\ds\int_{\mathbb{R}^{N}}V_0|\varphi_{\sigma,y}|^2\,dx &=&
\ds\int_{B_{\rho}(0)}V_{0}(x)|\varphi_{\sigma,y}|^2\,dx
+\ds\int_{\mathbb{R}^{N}\setminus
  B_{\rho}(0)}V_{0}(x)|\varphi_{\sigma,y}|^2\,dx
\\
&\leq& 
|V_{0}|_{L^{N/2}(B_{\rho}(0))}|\varphi_{\sigma,y}|^{2}_{L^{2^{*}}(B_{\rho}(0))}
+|V_{0}|_{L^{N/2}(\mathbb{R}^{N}\setminus
  B_{\rho}(0))}|\varphi_{\sigma,y}|^{2}_{L^{2^{*}}(\mathbb{R}^{N}\setminus
  B_{\rho}(0))}
\\
&\leq& 
|V_{0}|_{N/2}
\ds\sup_{y\in\mathbb{R}^{N}}|\varphi_{\sigma,y}|^{2}_{L^{2^{*}}(B_{\rho}(0))}
+|\varphi|^2_{2^*} |V_{0}|_{L^{N/2}(\mathbb{R}^{N}\setminus B_{\rho}(0))},
\qquad \forall \rho,\sigma>0.
\end{eqnarray*}
Using the fact that
$$
\ds\lim_{\sigma\rightarrow+\infty}\sup_{y\in\mathbb{R}^N}
|\varphi_{\sigma,y}|_{L^{2^{*}}(B_{\rho}(0))}=0, 
$$
we get
$$
\ds\lim_{\sigma\rightarrow+\infty}\sup\left\{
\ds\int_{\mathbb{R}^{N}}V_{0}(x)|\varphi_{\sigma,y}|^2\,dx;\
  \ y\in \mathbb{R}^{N}\ds\right\}\leq
|\varphi|^2_{2^*}\cdot |V_{0}|_{L^{N/2}(\mathbb{R}^{N}\setminus B_{\rho}(0))}. 
$$
Passing the limit of $\rho\rightarrow+\infty$ in the last inequality, we obtain $(b)$.

To prove $(c)$, we will assume by contradiction that there are sequences $(y_n)\subset\mathbb{R}^{N}$ and a sequence $(\sigma_n) \subset (0,+\infty)$ such that
\begin{equation}\label{L3-P2}
\ds\lim_{n\rightarrow+\infty}\int_{\mathbb{R}^{N}}V_{0}(x)\varphi^{2}_{\sigma_n,y_n}dx=L>0 \ \ {\rm and} \ \ |y_n|\rightarrow +\infty.
\end{equation}
From $(a)$ and $(b)$, we can suppose that
$$
\ds\lim_{n\rightarrow+\infty}\sigma_n=\overline{\sigma}>0.
$$
Using the hypotheses that $|y_n|\rightarrow+\infty$ and $V_{0}\in L^{N/2}(\mathbb{R}^{N})$, the Lebesgue's Theorem leads to 
$$
\ds\lim_{n\rightarrow+\infty}|V_{0}|_{L^{N/2}( B_{\sigma_{n}}(y_n))}=0,
$$
from where it follows that
$$
\ds\lim_{n\rightarrow+\infty}\int_{\mathbb{R}^{N}}V_{0}(x)\varphi^{2}_{\sigma_n,y_n}dx\leq \ds\lim_{n\rightarrow+\infty}|V_{0}|_{L^{N/2}(B_{\sigma_{n}}(y_n))}=0,
$$
which contradicts $(\ref{L3-P2})$. Therefore $(c)$ occurs.
{\fim}


\begin{lema}\label{LEMA4}
The following relations hold:
\begin{equation}
     \begin{cases}
     (a)\ \ \ds\lim_{\sigma\rightarrow 0}\sup\ds\left\{\gamma(\varphi_{\sigma,y});\ \ y\in\mathbb{R}^{N}\ds\right\}=0;\\
     (b)\ \ \ds\lim_{\sigma\rightarrow+\infty}\inf\ds\left\{\gamma(\varphi_{\sigma,y});\ \ y\in\mathbb{R}^{N},\ \ |y|\leq r\ds\right\}=1,\ \ \forall r>0;\\
     (c)\ \ (\beta(\varphi_{\sigma,y})|y)_{\mathbb{R}^{N}}>0;\ \ \forall y\in\mathbb{R}^{N}\setminus\{0\},\ \ \forall \sigma>0. 
     \end{cases}
\end{equation}   
Here $(x|y)$ denotes the usual inner product in $\mathbb{R}^N$ of the vectors $x,y \in \mathbb{R}^N$.  
\end{lema}
\noindent {\bf Proof.} Let $y\in\mathbb{R}^{N}$ be chosen arbitrarily. For any $\sigma>0$, using the fact that $\varphi_{\sigma,y} \in \mathcal{M}$ and the definitions of $\beta, \gamma$, we find
\begin{eqnarray}\label{LEMA4.P1}
0 &\leq&
\gamma(\varphi_{\sigma,y})=\frac{1}{|\varphi_{\sigma,y}|_{2^*}^{2^*}}
\ds\int_{\mathbb{R}^{N}}\ds\left|\frac{x}{1+|x|}-\beta(\varphi_{\sigma,y})
\ds\right||\varphi_{\sigma,y}|^{2^{*}}dx 
\nonumber\\
&\leq& 
\frac{1}{|\varphi_{\sigma,y}|_{2^*}^{2^*}}
\ds\int_{ B_{\sigma}(y)}\ds\left|\frac{x}{1+|x|}
-\frac{y}{1+|y|}\right||\varphi_{\sigma,y}|^{2^{*}}\,dx
+\ds\left|\ds\frac{y}{1+|y|}-\beta(\varphi_{\sigma,y})\ds\right|.
\end{eqnarray}
Now, from (\ref{PROPRIEDADE DA VARPHI}) and  (\ref{TRANSLACAO}), 
    \begin{eqnarray}\label{LEMA4.P2}
    \ds\left|\frac{y}{1+|y|}-\beta(\varphi_{\sigma,y})\ds\right| 
&=& 
\ds\frac{1}{|\varphi_{\sigma,y}|_{2^*}^{2^*}}\left|\int_{\mathbb{R}^{N}}
\ds\left(\ds\frac{y}{1+|y|}-\ds\frac{x}{1+|x|}
\ds\right)|\varphi_{\sigma,y}|^{2^{*}}dx\ds\right|
\nonumber\\
    &\leq& 
\ds\frac{1}{|\varphi_{\sigma,y}|_{2^*}^{2^*}}\int_{
  B_{\sigma}(y)}\ds\left|\ds\frac{y}{1+|y|}-
\ds\frac{x}{1+|x|}\ds\right||\varphi_{\sigma,y}|^{2^{*}}dx. 
\end{eqnarray}
Combining $(\ref{LEMA4.P1})$ with $(\ref{LEMA4.P2})$ we derive that
\begin{eqnarray*}
    0 \leq \gamma(\varphi_{\sigma,y}) &\leq&
    \frac{2}{|\varphi_{\sigma,y}|_{2^*}^{2^*}}\ds\int_{
      B_{\sigma}(y)}\ds\left|\frac{x}{1+|x|}-\frac{y}{1+|y|}
    \ds\right||\varphi_{\sigma,y}|^{2^{*}}dx
    \leq 2\sigma. 
\end{eqnarray*}
Hence
$$
0\leq \ds\sup\left\{\gamma(\varphi_{\sigma,y}); \ \
  y\in\mathbb{R}^{N}\right\}\leq 2\sigma, 
$$
which gives $(a)$ letting $\sigma\rightarrow 0$.

To prove $(b)$, let us first show that for all $y\in\mathbb{R}^{N}$,
\begin{equation}\label{LEMA4.P3}
    \ds\lim_{\sigma\rightarrow+\infty}\sup_{|y|\leq
      r}\left|\beta(\varphi_{\sigma,y})\right|=0. 
\end{equation}
Since $\varphi_{\sigma,0}$ is a symmetric function, we have
$\beta(\varphi_{\sigma,0})=0$. 
This combined with the limit below 
$$
\lim_{\sigma \to +\infty}\sup_{|y|\leq r}|\varphi_{1,y/\sigma} - \varphi_{1,0}|_{2^*}=0, 
$$
and the definition of $\beta$ gives (\ref{LEMA4.P3}).

 Now, fix $r>0$ arbitrarily and let $y\in\mathbb{R}^{N}$ such that
 $|y|\leq r$. 
For any $\sigma>0$, we see that
\begin{eqnarray*}
\gamma(\varphi_{\sigma,y}) 
&=& 
\frac{1}{|\varphi_{\sigma,y}|_{2^*}^{2^*}}\ds\int_{\mathbb{R}^{N}}\ds\left|
\ds\frac{x}{1+|x|}-\beta(\varphi_{\sigma,y})\ds\right||\varphi_{\sigma,y}|^{2^{*}}dx
\\
&\leq&
\frac{1}{|\varphi_{\sigma,y}|_{2^*}^{2^*}}\ds\int_{\mathbb{R}^{N}}
\ds\frac{|x|}{1+|x|}|\varphi_{\sigma,y}|^{2^{*}}dx+|\beta(\varphi_{\sigma,y})|
\\
&\leq& 1+|\beta(\varphi_{\sigma,y})|,
\end{eqnarray*}
which together with $(\ref{LEMA4.P3})$ leads us to
\begin{equation}\label{LEMA4.P4}
    \ds\limsup_{\sigma\rightarrow+\infty}\ds\left[\ds\inf\ds\left\{\gamma(\varphi_{\sigma,y}); \ \ y\in\mathbb{R}^{N}, \ \ |y|\leq r\ds\right\}\ds\right]\leq 1.
\end{equation}
If 
\begin{equation*}
\ds\limsup_{\sigma\rightarrow+\infty}\ds\left[\ds\inf\ds\left\{\gamma(\varphi_{\sigma,y}); \ \ y\in\mathbb{R}^{N}, \ \ |y|\leq r\ds\right\}\ds\right]< 1,
\end{equation*}
there are sequences $(\sigma_n)\subset (0,+\infty)$ and $(y_n)\subset\mathbb{R}^{N}$ such that $\sigma_n\rightarrow +\infty$, $|y_n|\leq r$ and
\begin{eqnarray}\label{LEMA4.P5}
\ds\lim_{n\rightarrow +\infty}\gamma(\varphi_{\sigma_{n},y_{n}})< 1.
\end{eqnarray}
On the other hand, considering $(\ref{LEMA4.P3})$, for all $\rho>0$ we deduce that
\begin{eqnarray*}
\gamma(\varphi_{\sigma_n,y_n}) &=& \frac{1}{|\varphi_{\sigma_n,y_n}|_{2^*}^{2^*}}\ds\int_{\mathbb{R}^{N}}\ds\left|\ds\frac{x}{1+|x|}-\beta(\varphi_{\sigma_n,y_n})\ds\right||\varphi_{\sigma_n,y_n}|^{2^{*}}dx\\
&\geq& \frac{1}{|\varphi_{\sigma_n,y_n}|_{2^*}^{2^*}} \ds\int_{\mathbb{R}^{N}}\ds\frac{|x|}{1+|x|}|\varphi_{\sigma_n,y_n}|^{2^{*}}dx-|\beta(\varphi_{\sigma_n,y_n})|\\
&\geq& \frac{1}{|\varphi_{\sigma_n,y_n}|_{2^*}^{2^*}} \ds\int_{\mathbb{R}^{N}\setminus B_{\rho}(0)}\ds\frac{|x|}{1+|x|}|\varphi_{\sigma_n,y_n}|^{2^{*}}dx-o_{n}(1)\\
&\geq& \ds\frac{\rho}{1+\rho}\frac{1}{|\varphi_{\sigma_n,y_n}|_{2^*}^{2^*}}\ds\int_{\mathbb{R}^{N}\setminus B_{\rho}(0)}|\varphi_{\sigma_n,y_n}|^{2^{*}}dx-o_{n}(1)\\
&\geq& \ds\frac{\rho}{1+\rho}\frac{1}{|\varphi_{1,0}|_{2^*}^{2^*}}\ds\int_{\mathbb{R}^{N}\setminus B_{\frac{\rho}{\sigma_n}}(-y_n/\sigma_n)}|\varphi_{1,0}|^{2^{*}}dx-o_{n}(1),
\end{eqnarray*}
hence
$$
\ds\lim_{n\rightarrow+\infty}\gamma(\varphi_{\sigma_n,y_n})\geq \ds\frac{\rho}{1+\rho}, \ \ \forall \rho>0.
$$
From this, since $\rho>0$ is arbitrarily,
$$
\ds\lim_{n\rightarrow+\infty}\gamma(\varphi_{\sigma_n,y_n})\geq 1,
$$
which contradicts $(\ref{LEMA4.P5})$. Thus, the equality in
$(\ref{LEMA4.P4})$ holds and the proof of $(b)$ is finished. 

Now, we will prove $(c)$. 
We note that if $0\notin B_{\sigma}(y)$, we have $(x|y)>0$ $\forall
x\in B_\sigma(y)$ and thus
$$
(\beta(\varphi_{\sigma,y})|y)=
\ds\int_{\mathbb{R}^{N}}\ds\frac{(x|y)}{1+|x|}|\varphi_{\sigma,y}|^{2^{*}}dx>0. 
$$
If $0\in B_{\sigma}(y)$, for each $x\in B_{\sigma}(y)$ such that
$(x|y)<0$, the point  $-x$  belongs to $B_{\sigma}(y)$, $(-x|y)>0$
and $\varphi_{\sigma,y}(-x)>\varphi_{\sigma,y}(x)$, 
which is enough to prove that $(\beta(\varphi_{\sigma,y})|y)>0$, as desired.
{\fim}

\begin{corolario}\label{CO-1}
There exist $r>0$ and $0<\sigma_1<\frac{1}{2}<\sigma_2$ such that
\begin{enumerate}

    \item [($a$)] $\gamma( {\varphi}_{\sigma_{1},y})<\ds\frac{1}{2},
      \quad \forall  y\in \mathbb{R}^{N}$, 
    
    \item [($b$)] $\gamma( {\varphi}_{\sigma_{2},y})>\ds\frac{1}{2},
      \quad \forall y\in \mathbb{R}^{N}$ with $|y|\le r$, 

\end{enumerate}
and 
\beq
\label{1825}
\ds\sup\left\{J_0( {\varphi}_{\sigma_{1},y})); \ \
  (\sigma,y)\in\partial \mathcal{H}\right\}<\overline{c}, 
\eeq
where
\begin{equation}
\label{K}
\mathcal{H}=[\sigma_1,\sigma_2]\times B_r(0).
\end{equation}
\end{corolario}
 
Points $(a)$ and $(b)$ follow from points   $(a)$ and $(b)$ of Lemma
$\ref{LEMA4}$, respectively, while    
(\ref{1825}) is a consequence of (\ref{1826}) and Lemma \ref{LEMA3}.
 
\begin{lema}\label{LEMA5}
Let $\mathcal{H}$ be the set defined in $(\ref{K})$.
Then, there exist $(\tilde{\sigma},\tilde{y})\in \partial\mathcal{H}$
and $(\overline{\sigma},\overline{y})\in\mbox{int}\,\mathcal{H}$
satisfying 
\beq
\label{1557}
\beta(\varphi_{\tilde{\sigma},\tilde{y}})=0 \ \ {\rm and} \ \
\gamma(\varphi_{\tilde{\sigma},\tilde{y}})>\ds\frac{1}{2} 
\eeq
and
\beq
\label{1558}
\beta(\varphi_{\overline{\sigma},\overline{y}})=0 \ \ {\rm and} \ \
\gamma(\varphi_{\overline{\sigma},\overline{y}})=\ds\frac{1}{2}.  
\eeq
\end{lema}
\noindent {\bf Proof.}\, 
First of all, note that by the symmetry of $\varphi$, we have
$\beta({\varphi}_{\sigma,0})=0$, $\forall\sigma>0$.  
Then $ (\tilde{\sigma},\tilde{y}):=(\sigma_2,0)$ verifies
(\ref{1557}), by Corollary \ref{CO-1}. 

In order to get (\ref{1558}) it is sufficient to consider also that
$\sigma\mapsto \gamma(\varphi_{\sigma,0})$ is a continuous map and
that $\sigma\mapsto \gamma(\varphi_{\sigma_1,0})<\frac{1}{2}$ while
$\sigma\mapsto \gamma(\varphi_{\sigma_2,0})>\frac{1}{2}$ by Corollary
\ref{CO-1}. 

{\fim}

\begin{lema}
\label{omotopia}
Let $g: \mathcal{H}\rightarrow \mathbb{R}\times
\mathbb{R}^{N}$ the function defined by
\beq
\label{1104}
g(\sigma,y)=(\gamma(\varphi_{\sigma,y}),\beta(\varphi_{\sigma,y})).
\eeq
Then,  
$$
deg(g,int(\mathcal{H}),(0,1/2)))=1.
$$
\end{lema}
\noindent {\bf Proof.}\, 
Let us consider the homotopy $G:[0,1]\times\partial
\mathcal{H}\rightarrow \mathbb{R}\times\mathbb{R}^{N}$ given by 
\beq
\label{1524}
G(t,\sigma,y)=(1-t)(\sigma,y)+
t(\gamma(\varphi_{\sigma,y}),\beta(\varphi_{\sigma,y})). 
\eeq
We remark $G$ is continuous and that
$$
G(0,\sigma,y)=(\sigma,y)
$$
and
$$
G(1,\sigma,y)=(\gamma(\varphi_{\sigma,y}),\beta(\varphi_{\sigma,y}))=g(\sigma,y).
$$
So, it remains to show that
\begin{equation}\label{GG}
\ds\left(\ds\frac{1}{2},0\ds\right)\notin G(t,\partial \mathcal{H}) \
\ \forall t\in [0,1] 
\end{equation}
or, equivalently,
$$
G(t,\sigma,y)\neq \ds\left(\ds\frac{1}{2},0\ds\right) \ \ \forall
(\sigma,y)\in\partial \mathcal{H} \ \ {\rm and} \ \ \forall t\in
[0,1]. 
$$
In fact, set $\partial \mathcal{H}=K_1\cup K_2\cup K_3$ with
\begin{equation*}
\begin{cases}    
K_1=\ds\left\{(\sigma,y); \ \ |y|\leq r, \ \ \sigma=\sigma_1\ds\right\},\\
K_2=\ds\left\{(\sigma,y); \ \ |y|\leq r, \ \ \sigma=\sigma_2\ds\right\},\\
K_3=\ds\left\{(\sigma,y); \ \ |y|=r, \ \ \sigma\in
  [\sigma_1,\sigma_2]\ds\right\}. 
\end{cases}
\end{equation*}
If $(\sigma,y)\in K_1$, then $\sigma=\sigma_1$ and by the Corollary
$\ref{CO-1}$ $(a)$ 
$$
(1-t)\sigma_1+t\gamma(\varphi_{\sigma_1,y})
<(1-t)\ds\frac{1}{2}+t\ds\frac{1}{2}
=\ds\frac{1}{2},
\ \ \forall t\in [0,1]. 
$$
Analogously, if $(\sigma,y)\in K_2$, then $\sigma=\sigma_2$ and again
by the Corollary $\ref{CO-1}$ $(b)$ 
$$
(1-t)\sigma_2+t\gamma(\varphi_{\sigma_2,y})>
(1-t)\ds\frac{1}{2}+t\ds\frac{1}{2}=\ds\frac{1}{2}, \ \ \forall t\in
[0,1]. 
$$
If $(\sigma,y)\in K_3$, then $|y|=r$ and
$0<\sigma_1\leq\sigma\leq\sigma_2$, so using Lemma $\ref{LEMA4}$
$(c)$, we obtain 
$$
\ds\left((1-t)y+t\beta(\gamma_{\sigma,y})|y\ds\right)
=(1-t)|y|^{2}+t(\beta(\varphi_{\sigma,y})|y)>0\qquad
\forall t\in[0,1]. 
$$
Now, the results follows by employing the proprieties of the Brouwer's
Topological degree.   
{\fim}


\begin{lema}\label{LEMA5bis}
Let $\mathcal{H}$ be the set defined in $(\ref{K})$, and assume that
$(V_3)$ holds, 
 then
$$
L=\max\{J_0({\varphi}_{\sigma,y}): (\sigma,y)\in \mathcal{H}\}
< 2^{\frac{4-\mu}{2N-\mu}}S_{H,L}.
$$
\end{lema}
\noindent {\bf Proof.} Using (\ref{familiadecatoff}), we have for all
$(\sigma,y)$ that 
\begin{eqnarray*}
J_0(  {\varphi}_{\sigma,y})&=&  \left[ \|\varphi_{\sigma,y}\|^{2}
  +\displaystyle\int_{\mathbb{R}^{N}}V_{0}(x)
  {\varphi}_{\sigma,y}^{2} dx \right] 
\leq \|\varphi\|^{2} +|V_{0}|_{L^{N/2}}|\varphi|_{2^*}^{2} \\
& \leq &  
\|\varphi\|^{2} +\frac{|V_{0}|_{L^{N/2}}}{S}\|\varphi\|^{2} \leq
\left[ 1+\frac{|V_{0}|_{L^{N/2}}}{S} \right]\Sigma \leq \left[
  1+\frac{|V_{0}|_{L^{N/2}}}{S} \right]\overline{c}. 
\end{eqnarray*}
The result follows by $(V_{3})$, (\ref{1930}) and (\ref{cbarra}). 

{\fim}



\section{Proof of Theorems}

Finally, with the help of the previous lemmas we are ready to prove
our main results. 
For $c\in\R$, let us fix the set 
$$
J^{c}_\lambda=\ds\left\{u\in\mathcal{M}:\ \ J_\lambda(u)\leq c\ds\right\}.
$$

\noindent {\bf Proof of Theorem $\ref{T1}$.} 
Combining the definition of $\overline{c}$ in (\ref{cbarra}), and
Lemma $\ref{LEMA5}$, we have 
$$
S_{H,L}< \overline{c}< \mathcal{B}_{V_0}\leq
J_{0}(\varphi_{\bar{\sigma},\bar{y}})
\leq L< 2^{\frac{4-\mu}{2N-\mu}}S_{H,L}.
$$

We will prove that functional $J_{0}$ constrained to $\mathcal{M}$ has
a critical level in the interval
$(S_{H,L},$ $2^{\frac{4-\mu}{2N-\mu}}S_{H,L})$.  
Suppose, by contradiction, that is not true. 
From Corollary \ref{LLL-2}, $J_{0}$ satisfies the Palais-Smale
condition in interval $(S_{H,L},2^{\frac{4-\mu}{2N-\mu}}S_{H,L})$.  
Thus, using a variant of the Deformation Lemma (see \cite{Struwe}) we
can find a  $\delta>0$ such that $\mathcal{B}_{V_{0}}-\delta>
{\bar{c}}$, $L+\delta<2^{\frac{4-\mu}{2N-\mu}}S_{H,L}$ and a continuous map
$\eta:J_{0}^{L+\delta}\rightarrow J_{0}^{\mathcal{B}_{V_0}-\delta}$
such that 
$$
\eta(u)=u, \ \ \forall u\in J_{0}^{\mathcal{B}_{V_0}-\delta}.
$$
Then,  the map $\eta(\varphi_{\sigma,y})$, $(\sigma,y)\in\mathcal{H}$,
is well defined and  we remark that 
$$
 J_{0}(\eta(\varphi_{\sigma,y}))<\mathcal{B}_{V_0}-\delta, \quad
 \forall (y,\sigma)\in \mathcal{H},  
$$
which implies
\beq
\label{1105}
\Theta(\varphi_{\sigma,y})
:=(\gamma(\eta(\varphi_{\sigma,y})),\beta(\eta(\varphi_{\sigma,y})))\neq
\left(\frac 12,0\right). 
\eeq
On the other hand,  by Corollary \ref{CO-1},
\begin{eqnarray}\label{pardarcertomaisumavez}
 J_{0}(\varphi_{\sigma,y})<\overline{c}<\mathcal{B}_{V_0}-\delta,
 \quad \forall (y,\sigma)\in \partial\mathcal{H},   
\end{eqnarray}
which implies $\eta(\varphi_{\sigma,y})=\varphi_{\sigma,y}$, from which
$$
\Theta(\varphi_{\sigma,y})=
g(\sigma,y)=(\gamma(\varphi_{\sigma,y}),\beta(\varphi_{\sigma,y})),
\quad \forall (y,\sigma)\in \partial\mathcal{H}, 
$$
where $g$ is the map introduced in (\ref{1104}).

Therefore, by the homotopy invariance of topological degree, taking
into account Lemma \ref{omotopia}, we deduce 
$$
1=d(g, \mathcal{H},( 1/2,0))=d(\Theta, \mathcal{H},(1/2,0)),
$$
that implies the existence of $(\sigma,y)$ such that
$\Theta(\sigma,y)=( 1/2,0)$, contradicting
(\ref{1105}). 
Therefore, the functional $J_{0}$ constrained on $\mathcal{M}$ has at
least one critical point $u\in\mathcal{M}$ such that
$\overline{c}<J(u)<2^{\frac{4-\mu}{2N-\mu}}S_{H,L}$. 
Moreover, by Lemma $\ref{Nodal}$, we also have $u>0$, finishing
the proof.

\noindent {\bf Proof of Theorem $\ref{T2}$.} 

The proof of this results follows as in \cite{CM}, however for the
reader's convenience we will write its proof.  

\begin{lema}\label{LEMAfinal}
Let $\mathcal{H}$ be the set defined in $(\ref{K})$. Then, there
exists $\bar\lambda>0$ such that, for all $0<\lambda<\bar\lambda$, we
have 
$$
l_\lambda:=\max\{J_\lambda({\varphi}_{\sigma,y}):
(\sigma,y)\in \partial\mathcal{H}\}< \overline{c}. 
$$
Moreover, $\gamma(\varphi_{\sigma_1,y})<\frac 12$,
$\gamma(\varphi_{\sigma_2,y})>\frac 12$, for all $y \in
\mathbb{R}^{N}$, $|y|< r$. 
\end{lema}
\noindent {\bf Proof.} Note that 
$$
\lambda \displaystyle\int_{\mathbb{R}^{N}} |{\varphi}_{\sigma,y}|^{2}
dx=\lambda \sigma^{2} \displaystyle\int_{B_{1}(0)}|\varphi|^2 dx, 
$$
and so, 
$$
J_\lambda(  {\varphi}_{\sigma,y})=J_0(  {\varphi}_{\sigma,y})+\lambda
\sigma^{2}   \int_{B_{1}(0)}|\varphi|^2 dx. 
$$
Now, the result follows from (\ref{1825}). 
{\fim}

Combining the definition of $\overline{c}$ and Lemma 
$\ref{LEMAfinal}$, for every $\lambda\in(0,\bar\lambda)$ we have
\beq
\label{73}
S_{H,L}<\Upsilon_\lambda \leq
J_{\lambda}(\varphi_{\tilde{\sigma},\tilde{y}})\leq 
l_\lambda <\overline{c}< \mathcal{B}_{V_0}, 
\eeq
where $\Upsilon_\lambda$ has been defined in (\ref{1145}) and
$(\tilde\sigma,\tilde y)$ has been introduced in Lemma \ref{LEMA5}.

We will prove that functional $J_{\lambda}$ constrained to
$\mathcal{M}$ has a critical level in the interval $(\Upsilon_\lambda
,l_\lambda )$.  
Suppose, by contradiction, that is not true. 
From Corollary \ref{LLL-2}, $J_{\lambda}$ satisfies the Palais-Smale
condition in interval $(\Upsilon_\lambda,l_\lambda)$. 
Thus, using a variant of the Deformation Lemma (see \cite{Struwe}) we
can find a  
 positive number $\delta_1>0$ such that   $\Upsilon_\lambda-\delta_{1}> S_{H,L}$,
$l_\lambda +\delta_{1}<\overline{c}$ and a continuous function
$$
\eta:[0,1]\times J_\lambda^{l_\lambda +\delta_{1}}\lo J_\lambda^{l_\lambda+\delta_1}
$$
such that
\begin{eqnarray}
\nonumber
& \eta(0,u)=u  & \hspace{1cm}\forall u\in J_\lambda^{l_\lambda+\delta_1}
\\
\nonumber
& \eta(s,u)=u &\hspace{1cm}\forall u\in
J_\lambda^{\Upsilon_\lambda-\delta_{1}}, \ \forall s\in[0,1] 
\\
\label{44.1}
&J_\lambda( \eta\, (s,u))\le J_\lambda(u) &\hspace{1cm} \forall u\in
J_\lambda^{\Upsilon_\lambda+\delta_{1}}, \  \forall s\in[0,1] 
\\
\label{44.2}
& \eta(1,J_\lambda^{l_\lambda +\delta_{1}} )\subset J_\lambda^{\Upsilon_\lambda-\delta_{1}}. &
\end{eqnarray}
Therefore, definition of $ l_\lambda$ and (\ref{44.2}) give 
\beq
\label{44.3}
(\sigma,y)\in\partial\cH\ \Rightarrow\
J_\lambda(\varphi_{\sigma,y})\le l_\lambda\ \Rightarrow\
J_\lambda(\eta(1,\varphi_{\sigma,y})\le\Upsilon-\delta_1.
\eeq
Let us consider $\forall s\in [0,1]$, $\forall (\sigma,y)\in\partial \cH$
$$
\Gamma(\sigma,y,s)\ =\ \left\{
\begin{array}{lc}
G(\sigma,y,2s), & s\in[0,1/2] \\
\left(\gamma\circ\eta\, 
  (2s-1),\varphi_{\sigma,y}),\beta\circ\eta\, (2s-1,\varphi_{\sigma,y})\right)
& s\in[1/2,1],
\end{array}
\right.
$$
where $G$ is the map defined in (\ref{1524}). 
As already shown in Lemma \ref{omotopia},
\beq
\label{111}
\forall s\in[0,1/2],\ \forall (\sigma,y)\in\partial \cH,\
G(\sigma,y,s)\neq \left( \frac{1}{2},0\right). 
\eeq
Furthermore, from (\ref{73}) and (\ref{44.1}), we deduce 
$\forall
s\in[1/2,1]$ $\forall (\sigma,y)\in\partial\cH$
$$
J_\lambda(\eta(2s-1,\varphi_{\sigma,y}))\le
J_\lambda(\varphi_{\sigma,y})\le  l_\lambda <\bar
c<\cB_{V_0}\le\cB_{V_\lambda},\quad \forall \lambda>0, 
$$
which gives
\beq
\label{222}
\forall s\in[1/2,0],\ \forall (\sigma,y)\in\partial \cH,\ G(\sigma,y,s) \neq \left( \frac{1}{2},0\right). 
\eeq
By (\ref{111}), (\ref{222}) and the continuity of $\Gamma$, we obtain the existence of $(\check \sigma,\check
y)\in\partial\cH$ such that 
$$
\gamma\circ\eta\, (1,\varphi_{\check \sigma,\check y})\ge \frac{1}{2},\quad 
\beta\circ\eta\, (1,\varphi_{\check \sigma,\check y})=0.
$$
Then
$$
J_\lambda( \eta(1,\varphi_{\check \sigma,\check y}))\ge \Upsilon_\lambda,
$$
which contradicts (\ref{44.3}).
Therefore, the functional $J_{\lambda}$
constrained on $\mathcal{M}$ has at least one critical point
$u_l\in\mathcal{M}$, such 
that $S_{H,L}<J(u_l)<\overline{c}$, $\forall \lambda\in(0,\bar \lambda)$. 
Moreover, by the Lemma $\ref{Nodal}$, we deduce $u_l>0$, concluding
the proof of the first part of the theorem.  
 
%
%
%

\vspace{1mm}

Now, let us suppose that $(V_3)$ holds.
Then the existence of an high energy positive solution, i.e. of a
critical point for $J_\lambda$ constrained on $\mathcal{M}$ such that
$\bar c\le J_\lambda(u_h)< 2 S_{H,L}$, can be proved for small
$\lambda$ arguing exactly as in the proof 
of Theorem \ref{T1}, taking into account that    
$$
\lim_{\lambda\to 0} 
\sup\{J_\lambda({\varphi}_{\sigma,y}): (\sigma,y)\in \mathcal{H}\}
=\sup\{J_0({\varphi}_{\sigma,y}): (\sigma,y)\in \mathcal{H}\}
$$ 
and  
$$
\lim_{\lambda\to 0} 
\sup\{J_\lambda({\varphi}_{\sigma,y}): (\sigma,y)\in \partial \mathcal{H}\}
=\sup\{J_0({\varphi}_{\sigma,y}): (\sigma,y)\in  \partial
\mathcal{H}\}.
$$

\vspace{0.3cm}
\noindent {\bf Acknowledgment:} 
C.O. Alves is partially supported by CNPq/Brazil 304804/2017-7,  
G. M. Figueiredo is supported by CNPq and FAPDF,
R. Molle is supported by the MIUR Excellence Department Project CUP
E83C18000100006 (Roma Tor Vergata University) and by the INdAM-GNAMPA
group.

\end{document}